\documentclass{article}

\usepackage[margin=2.5cm]{geometry}

\usepackage[utf8]{inputenc} 
\usepackage[title]{appendix}

\usepackage[cmex10]{amsmath}
\usepackage{amssymb,amsfonts}
\usepackage{bbm} 
\usepackage{amsthm} 

\newtheorem{theorem}{Theorem}
\newtheorem{proposition}{Proposition}
\newtheorem{definition}{Definition}
\newtheorem{lemma}{Lemma}
\newtheorem{corollary}{Corollary}
\newtheorem{assumption}{Assumption}

\DeclareUnicodeCharacter{2212}{-}

\usepackage{subfig}

\usepackage{hyperref}
\hypersetup{
  colorlinks   = true, 
  urlcolor     = blue, 
  linkcolor    = blue, 
  citecolor   = blue 
}
\usepackage{tikz}
\usepackage{tikzscale}
\usetikzlibrary{patterns}
\usetikzlibrary{datavisualization.formats.functions}
\usetikzlibrary{decorations.pathreplacing,calligraphy}
\tikzstyle{solve}=[shape=rectangle, rounded corners, draw, align=center, fill=white]
\tikzstyle{condition}=[draw=black, align=center, fill=white]

\newcommand{\uptxt}[1]{^{\mathrm{#1}}}

\newif\ifarxiv
\arxivtrue

\usepackage{authblk}
\usepackage[round]{natbib}

\usepackage{setspace}
\title{Sufficient A Priori Conditions for the Linear Relaxation of the Energy Storage Scheduling Problem}
\author[1]{Eléa Prat (emapr@dtu.dk)}
\author[1]{Richard Martin Lusby (rmlu@dtu.dk)}
\author[2,1,3]{Pierre Pinson (p.pinson@imperial.ac.uk)}
\affil[1]{Department of Technology, Management and Economics, Technical University of Denmark}
\affil[2]{Dyson School of Design Engineering, Imperial College London}
\affil[3]{Halfspace, Denmark} 

\begin{document}

\maketitle

\begin{abstract}
    When modeling energy storage systems, an essential question is how to account for the physical infeasibility of simultaneous charge and discharge. The use of complementarity constraints or of binary variables is common, but these formulations do not scale well. Alternatively, assumptions such as perfect efficiencies or positive prices are often used to justify the choice of a linear model. In this paper, we establish new a priori conditions that guarantee the existence of an optimal solution without simultaneous charge and discharge when solving the linear relaxation of the storage scheduling problem. They are based on the characteristics of the storage system, in particular, the duration of charge. They can be valid for negative prices and with inefficiencies, thereby enlarging the set of conditions for which the complementarity constraints can be relaxed. We prove mathematically the validity of these conditions and illustrate them with practical examples.
We also introduce a refined mixed-integer linear equivalent, in which the number of binary variables can be drastically reduced.
\end{abstract}

\section{Introduction}

The significant efforts to promote renewable energies bring many challenges, in particular, due to their intermittent nature.
Energy storage has established itself as a desirable solution to counter some of these effects.
In this context, accurately modeling such storage systems is highly relevant and many questions are still to be addressed.
One key challenge is to determine the length of the solving horizon for energy-storage problems, see \cite{Cruise2019Control}, while another is to enhance the efficiency of solving these problems \citep{Sioshansi2021Energy}.
Energy storage systems are featured in various relevant applications, such as the scheduling of a system with renewable energy production and storage \citep{Kim2011Optimal}, or of a standalone storage system participating in a market \citep{Cruise2019Control}.

A challenge when modeling energy storage systems is to account for the exclusivity of charging and discharging modes \citep{Chen2023Security}, as solutions with simultaneous charge and discharge are not physically feasible.
This requirement can be modeled with a complementarity constraint, imposing that the product of charged and discharged quantities equals zero. However, this renders the problem non-linear. 
Alternatively, this problem can be reformulated as a Mixed-Integer Linear Programming (MILP) problem, in which binary variables control charging and discharging modes. This option can still be computationally intractable when the number of binary variables is high. Moreover, it is non-convex, making it unsuitable for some applications.
Therefore, a common approach is to solve the linear relaxation of the problem, in which the complementarity constraints are removed.
In this context, a complete analysis of when such a relaxation returns a solution that does not exhibit simultaneous charge and discharge is lacking and would be extremely valuable for informing on the most appropriate mathematical formulation.

Apart from these approaches, other formulations and heuristics have been identified in the literature.
Several recent works propose the use of robust formulations, which are linear and ensure the exclusivity of charge and discharge modes, but at the expense of reducing the solution space -- see \cite{Pozo2022Linear, Zhang2022Pricing, Yildiran2023Robust, Lin2024ACFeasible}.
The solution obtained can be significantly conservative, e.g., when the charging and discharging efficiencies are low, as shown by \cite{Nazir2021Guaranteeing}.
Another option is to solve the linear relaxation of the problem, in which the complementarity constraints are removed and modify the solution to one without simultaneous charge and discharge -- see \cite{Almassalkhi2014Model, Tziovani2023Successive}.
\cite{Elsaadany2023Battery} compare different models for scheduling a storage system and show that the solutions obtained by relaxing and then correcting in case of simultaneous charge and discharge can also be far from optimal. 

Due to the convenience of the linear relaxation of the storage model, in recent years, progress has been made in identifying conditions under which the relaxation of complementarity constraints is \emph{exact}, i.e., that there exists a solution to this relaxation that does not exhibit simultaneous charge and discharge.
However, most work focusing on establishing such conditions delivers \emph{a posteriori conditions}, see e.g. \cite{Garifi2020Convex, Nazir2018Receding, Nazir2020Optimal, Li2015Sufficient, Li2018Extended, Duan2016Improved, Baldick2023Optimization, Wang2023Sufficient}. These are impractical \citep{Yildiran2023Robust}.

The only a priori conditions to date assume the perfect efficiency of the storage system \citep{Almassalkhi2014Model,Yildiran2023Robust} or the positivity of prices \citep{Haessig2021Convex}. Both assumptions are very restrictive in practice, as real storage systems usually suffer non-negligible losses when charging and discharging, and electricity market prices, such as day-ahead prices, can take negative values. In fact, the occurrence of strictly negative prices is increasing with the penetration of renewable energy sources. Moreover, it is essential to properly account for negative prices as they can significantly modify the optimal strategy \citep{Zhou2016Electricity}.

The energy storage system is an example of a capacitated inventory or warehouse problem \citep{Charnes1955Generalizations} with losses, complementarity constraints, and bounds on sales and purchases. 
Therefore, the development of energy storage also motivates advances in the inventory and warehouse literature.
\cite{Wolsey2018Convex} present a convex hull for the capacitated warehouse problem with complementarity constraints. However, sales and purchases are only limited by the stored quantity and the storage capacity, respectively.
\cite{Bansal2023Warehouse} consider time-dependent bounds on sale and purchase quantities, making it harder to retrieve general results regarding the possible relaxation of complementarity constraints.

In this paper, we consider the problem of optimally scheduling a storage system to maximize profit, addressing the question: Can we extend the conditions for which we know that the linear relaxation will return an optimal solution to the storage scheduling problem with complementarity constraints?

To answer this question, our first contribution is to formalize the results for positive prices and perfect efficiency with Proposition~\ref{prop1}.
A second contribution is the introduction of a refined MILP relaxation of the storage scheduling problem, which has significantly fewer binary variables.
Finally, our main contributions are the mathematical proofs of different sufficient conditions for which we can guarantee the (in)exactness of the linear relaxation \emph{a priori}, and even when considering storage inefficiencies and negative prices.
In Theorem \ref{theorem1}, we give a sufficient condition based on the longest sequence of strictly negative prices for which we know that simultaneous charge and discharge \emph{will be} part of any optimal solution. On the other hand, Theorem \ref{theorem2} introduces sufficient conditions under which the linear relaxation is exact for a special case, and Theorem~\ref{theorem3} is a generalization of this.
Though the stylized model used is motivated by energy storage, the results are valid for any other inventory problems with an equivalent setup.

The paper is organized as follows. In Section \ref{sec:model}, we present the main concepts and the state of the art.
Section \ref{sec:ideal} formalizes existing results and introduces a refined MILP relaxation of the problem. In Section~\ref{sec:new}, we demonstrate new, a priori, sufficient conditions for the (in)exact relaxation of complementarity constraints and illustrate them with examples. Section \ref{sec:discussion} describes how the different findings can be used and discusses some limitations. Section \ref{sec:ccl} concludes. Proofs of the results in the main text are given as appendices.

\section{Problem Set-Up} \label{sec:model}

We consider the problem of scheduling a price-taker energy storage system over a finite time interval~$\mathcal{T}=\{1,2,...,T\}$, where the duration of each time period is $\Delta t$.

\subsection{Representation of the Storage System} \label{sec:model_stg}

We use a simplified representation of the storage system \citep{Pozo2014Unit,Pozo2022Linear} as described below, and consider that the storage system can be continuously dispatched.
This representation has limitations, see \cite{Sioshansi2021Energy}, but provides a good representation of the charging and discharging dynamics, and for this reason, is largely used in the literature on storage system scheduling. 
\ifarxiv
\begin{figure}[h]
    \centering
    {\resizebox{.13\linewidth}{!}{\includegraphics{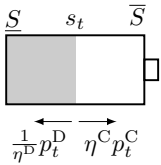}}}
    \caption{Schematic representation of the storage system and associated notations}
    \label{fig:stg}
\end{figure}
\fi

The decision variables are the storage system state of energy, $s_t$, and the power charged and discharged, $p_t\uptxt{C}$, and $p_t\uptxt{D}$. We denote by $\mathbf{x}$ the vector of all variables, over all the time periods.
The variables are represented in Figure \ref{fig:stg}. The figure also shows the lower and upper bounds on the storage state of energy, $\underline{S}\geq0$ and $\overline{S}>\underline{S}$, as well as the charging and discharging efficiencies, $\eta\uptxt{C} \in (0,1]$ and $\eta\uptxt{D} \in (0,1]$. The round-trip efficiency is $\eta = \eta\uptxt{C} \eta\uptxt{D}$. Due to self-discharge, the energy available at the beginning of a time period is a fraction $\rho$ of the energy left at the end of the previous time period, where $\rho \in (0,1]$. The state of energy at the end of each time period can then be expressed as
\begin{equation} \label{eq:update}
    s_t = \rho s_{t-1} +  \Delta t (\eta\uptxt{C} p_t\uptxt{C} - \frac{1}{\eta\uptxt{D}} p_t\uptxt{D}) \, , \, \forall t \in \mathcal{T} \setminus \{1\} \, .
\end{equation}

For the first time period, we need to consider the initial level $S\uptxt{init}$:
\begin{equation} \label{eq:init}
    s_1 = \rho S\uptxt{init} +  \Delta t (\eta\uptxt{C} p_1\uptxt{C} - \frac{1}{\eta\uptxt{D}} p_1\uptxt{D}) \, .
\end{equation}

The final state of energy is left free.
In Section \ref{sec:discussion}, we discuss the impact of constraining it.

The following constraints are needed to set the bounds. First, for the state of energy we have
\begin{equation}
    \label{eq:bounds_s} \underline{S} \leq s_t \leq \overline{S} \, , \, \forall t \in \mathcal{T} \,.
\end{equation}
For the quantity charged, with maximum charging rate, $\overline{P}\uptxt{C}>0$, the bounds are 
\begin{equation}
    \label{eq:bounds_pc} 0 \leq p_t\uptxt{C} \leq \overline{P}\uptxt{C} \, , \, \forall t \in \mathcal{T} \,.
\end{equation}
Finally, for the quantity discharged, with maximum discharging rate, $\overline{P}\uptxt{D}>0$, the bounds are 
\begin{equation}
    \label{eq:bounds_pd} 0 \leq p_t\uptxt{D} \leq \overline{P}\uptxt{D} \, , \, \forall t \in \mathcal{T} \,.
\end{equation}

In our derivations, we use the fact that the state of energy at any $t_2$ can be expressed using the state of energy at any previous time period $t_1 < t_2$, by merging constraints \eqref{eq:update} for all $t$ between $t_1+1$ and $t_2$. This gives:
\begin{equation} \label{eq:s_calc}
    s_{t_2} = \rho^{t_2-t_1} s_{t_1} + \Delta t\sum_{t=t_1+1}^{t_2} \rho^{t_2-t}( \eta\uptxt{C} p_t\uptxt{C} - \frac{1}{\eta\uptxt{D}} p_t\uptxt{D} ),
\end{equation}
with $s_0=S\uptxt{init}$, when $t_1=0$.

Finally, the following assumption ensures that the quantity lost by leakage can always be recovered by charging the maximum quantity.

\begin{assumption}\label{ass3}
    $(1-\rho) \overline{S} \leq \Delta t \, \eta\uptxt{C} \overline{P}\uptxt{C} $.
\end{assumption}

The main results of this paper are based on the characteristics of the storage system, and in particular, on how long it takes to fully charge the storage system. 
\begin{definition}[Duration and speed of charge / Duration and speed of discharge]
The \emph{speed of charge} of a storage system corresponds to the maximum charging rate, accounting for inefficiencies, $\eta\uptxt{C} \overline{P}\uptxt{C}$.
The \emph{duration of charge} of a storage system, $\Delta t \uptxt{C}$, is the minimum duration required to fully charge the storage system, charging from the minimum level and disregarding leakage:
\begin{equation*}
    \Delta t \uptxt{C} = \frac{\overline{S} - \underline{S}}{\eta\uptxt{C} \overline{P}\uptxt{C}}.
\end{equation*}
Similarly, the \emph{speed of discharge} is $({1}/{\eta\uptxt{D}}) \overline{P}\uptxt{D}$.
The \emph{duration of discharge}, $\Delta t \uptxt{D}$, is
\begin{equation*}
    \Delta t \uptxt{D} = \frac{\overline{S} - \underline{S}}{({1}/{\eta\uptxt{D}}) \overline{P}\uptxt{D}}.
\end{equation*}
\end{definition}

In the subsequent examples, we compare results for different durations of charge and discharge. Since we assume a fixed capacity for the storage systems in these examples, the parameters that vary are those associated with the speed of charge and discharge. We use the qualifications of ``fast'' and ``slow'' storage to identify an example where the storage system charges \emph{faster} than in the other example in which the storage system charge is \emph{slower}, and usually compared to a threshold value defined by one of our different results.

When scheduling energy storage systems, we must take into account the fact that solutions exhibiting simultaneous charge and discharge are not physically feasible and must be avoided.
The mathematical definition of simultaneous charge and discharge is as follows.
\begin{definition}[Simultaneous charge and discharge]
A solution $\mathbf{x}^*$ to a problem with storage -- modeled following \eqref{eq:update} to \eqref{eq:bounds_pd} -- exhibits simultaneous charge and discharge if $ \exists \, t \in \mathcal{T}$ such that $p_t\uptxt{C*} > 0$ and $p_t\uptxt{D*} > 0$.
\end{definition}


\subsection{Models}

The objective is to maximize the profit from arbitrage using the storage system. The energy price paid if charging, or received if discharging, is $C_t$. We further assume that the storage system considered is a small player in the market, acting as a price-taker.
There will always be enough demand or production to fulfill the desired schedule. The objective function to maximize is then:
\begin{equation} \label{eq:obj} 
    Z = \sum_{t \in \mathcal{T}} \Delta t \, C_t (p_t\uptxt{D}-p_t\uptxt{C}).
\end{equation}

We first present the model with complementarity constraints and its equivalent MILP formulation. We then introduce the linear relaxation in which the complementarity constraints are dropped.

\subsubsection{Non-Linear Model.}
The model with complementarity constraints is
\begin{subequations} \label{pb:opt_nl}
\begin{align}
    \label{eq:nl_obj} \max_{\mathbf{x}} \quad & Z \\
    \label{eq:nl_stg_mdl} \text{s.t.} \quad & \eqref{eq:update} - \eqref{eq:bounds_pd} \, , \\
    \label{eq:complementarity} & p_t\uptxt{C} p_t\uptxt{D} = 0\, ,& \forall t \in \mathcal{T}.
\end{align}
\end{subequations}
Constraint \eqref{eq:complementarity} is introduced to exclude solutions with simultaneous charge and discharge. 
We denote by $\mathcal{X}\uptxt{CC}$ and $\mathcal{X}\uptxt{CC*}$ the sets of feasible and optimal solutions, respectively, to \eqref{pb:opt_nl}.

\subsubsection{MILP Formulation.}
The complementarity constraints can be linearized using binary variables $u_t\uptxt{C}$ and $u_t\uptxt{D}$, to model the choice between charge ($u_t\uptxt{C}=1$) and discharge ($u_t\uptxt{D}=1$). The MILP model is then:
\begin{subequations}  \label{pb:opt_milp}
\begin{align}
    \label{eq:milp_obj} \max_{\mathbf{x}} \quad & Z  \\
    \label{eq:milp_stg_mdl} \text{s.t.} \quad & \eqref{eq:update} - \eqref{eq:bounds_s} \, , \\
    \label{eq:bounds_pc_bin} & 0 \leq p_t\uptxt{C} \leq u_t\uptxt{C} \overline{P}\uptxt{C} \, ,& \forall t \in \mathcal{T}, \\
    \label{eq:bounds_pd_bin} & 0 \leq p_t\uptxt{D} \leq u_t\uptxt{D} \overline{P}\uptxt{D} \, ,& \forall t \in \mathcal{T}, \\
    \label{eq:complementarity_bin} & u_t\uptxt{C} + u_t\uptxt{D} \leq 1 \, ,& \forall t \in \mathcal{T}, \\
    \label{eq:bin_def} & u_t\uptxt{C}, u_t\uptxt{D} \in \{ 0,1\} \, ,& \forall t \in \mathcal{T}.
\end{align}
\end{subequations}
Constraint \eqref{eq:complementarity_bin} prevents simultaneous charge and discharge. Constraints \eqref{eq:bounds_pc_bin} and \eqref{eq:bounds_pd_bin} are modifications of~\eqref{eq:bounds_pc} and \eqref{eq:bounds_pd} to force the variables to zero when the corresponding binary variable, introduced in \eqref{eq:bin_def}, is zero. 

\subsubsection{Linear Relaxation.}

The linear relaxation of \eqref{pb:opt_nl} is
\begin{subequations} \label{pb:opt}
\begin{align}
    \label{eq:lp_obj} \max_{\mathbf{x}} \quad & Z \\
    \label{eq:lp_stg_mdl} \text{s.t.} \quad & \eqref{eq:update} - \eqref{eq:bounds_pd} \, .
\end{align}
\end{subequations}
The sets of feasible and optimal solutions of \eqref{pb:opt} are $\mathcal{X}\uptxt{LP}$ and $\mathcal{X}\uptxt{LP*}$, respectively.


We want to derive a priori conditions under which there is a solution of \eqref{pb:opt} that is feasible, and thus optimal, for the original problem in \eqref{pb:opt_nl}. 

\begin{definition}[Exact and inexact relaxation]
The relaxation \eqref{pb:opt} is \emph{exact} if $\mathcal{X}\uptxt{LP*} \cap \mathcal{X}\uptxt{CC} \neq \varnothing$. Otherwise, the relaxation \eqref{pb:opt} is \emph{inexact}.

\end{definition}

The relaxation is exact if there exists a solution $\mathbf{x}^* \in \mathcal{X}\uptxt{LP*}$ that does not exhibit simultaneous charge and discharge. The relaxation is inexact if every solution $\mathbf{x}^* \in \mathcal{X}\uptxt{LP*}$ exhibits simultaneous charge and discharge.

\subsection{Related Work} \label{sec:related}

Many of the sufficient conditions established by earlier works are a posteriori conditions that can only be verified after solving the problem.
They are either based on the value of dual variables, in particular with the (strict) positivity of locational marginal prices (LMPs), or on the value of primal variables, checking that the storage system's bounds are not reached \citep{Garifi2020Convex, Castillo2013Profit, Duan2016Improved}.
Conditions based on the value of the dual variables are given by \cite{Nazir2018Receding} and \cite{Nazir2020Optimal}. 
For the economic dispatch problem, \cite{Li2015Sufficient} show that if the discharging prices are all higher than the charging prices and if the charging costs are strictly lower than the LMP at the storage location, the complementarity constraint can be relaxed. 
A milder condition is given by \cite{Duan2016Improved}, who highlight that the non-negativity of LMPs is not a requirement for the relaxation to be exact. Two more conditions based on LMPs are given by \cite{Li2018Extended}. 
New conditions on LMPs introduced by \cite{Wang2023Sufficient} are proven to be more general than six of the conditions in the previously cited works.
They also include a penalty in the objective function to discourage simultaneous charge and discharge, and show the impact of setting the penalty on the validity of the relaxation. A similar approach is followed by \cite{Han2025Regularized} and by \cite{Zarilli2018Energy}.
In several of these papers, it is argued that the LMPs can be forecast, but this would only give a likelihood that the condition will be valid.

Closer to our work, \cite{Baldick2023Optimization} derive results for a finite horizon model of a pumped storage hydroelectric system scheduling. The first concerns the suboptimality of simultaneously charging and discharging in time periods with strictly positive prices. However, the case of prices equal to zero and efficiencies equal to one is not covered. The second relates simultaneous charge and discharge to reaching the capacity bounds of the storage system. Here again, there is no insight into how this can be determined before solving the problem. Valid inequalities are introduced for the storage scheduling problem, similar to those by \cite{Chen2022Battery}, which show that the valid inequalities can remove some of the solutions with simultaneous charge and discharge, but they are not sufficient to completely avoid it.

Other works focus on a priori conditions, but these are very restrictive. Sufficient conditions for an optimal portfolio problem including storage are given by \cite{Yang2014Joint}, who show that if the cost function for charging and discharging is increasing and positive, there will not be simultaneous charge and discharge.
A paper by \cite{Almassalkhi2014Model} gives an expression for the difference in the state of energy between a model with complementarity constraints and a model where they are relaxed. They show that the two models match when the round-trip efficiency is one. 
Sufficient conditions are also given by \cite{Rahbar2015Realtime}, \cite{Wu2017Distributed} and \cite{Joshi2021Sufficient}. They are based on the fact that prices are strictly positive and the round-trip efficiency is strictly lower than one.
In \cite{Li2015StorageLike}, sufficient a priori conditions are identified, but these are only valid for the problem of load leveling, and not for the storage scheduling problem studied here.
In \cite{Lin2023Relaxing}, a condition based on the storage size is presented, namely that the storage system can be net charging the maximum quantity over the full horizon without violating the capacity. This is much more restrictive than the conditions on storage size that we introduce.

\subsection{Relevant Subsets}
The conditions introduced next rely on different subsets of time periods, which are defined here.

We identify $\mathcal{T}^-=\{t \in \mathcal{T} \,|\, C_t <0\}$, the subset of time periods with strictly negative prices, and $\mathcal{T}^+=\{t \in \mathcal{T} \,|\, C_t \geq 0\}$, the subset of time periods with positive prices, such that $\mathcal{T} = \mathcal{T}^- \cup \mathcal{T}^+$.
We also identify $\mathcal{T}^0=\{t \in \mathcal{T} \,|\, C_t =0\}$, the subset of time periods with a price of zero.
We denote by $\overline{\mathcal{T}}^-$ the compact subset of $\mathcal{T}^-$ with the highest cardinality, $\overline{n} = |\overline{\mathcal{T}}^-|$ .

$\mathcal{T}$ can be stated as an ordered collection of compact subsets of $\mathcal{T}^+$ and $\mathcal{T}^-$, which are identified as $\mathcal{T}^+_j$ and $\mathcal{T}^-_j$, with $j \in \mathcal{J}$, such that $\mathcal{T} = \bigcup_{j \in \mathcal{J}} ( \mathcal{T}^+_j \cup \mathcal{T}^-_j )$. It is therefore a representation of $\mathcal{T}$ as a collection of consecutive time periods with positive prices and consecutive time periods with strictly negative prices.
Note that if $C_1<0$, $1\in\mathcal{T}^-_1$ and $\mathcal{T}^+_1 = \varnothing$. 
We also introduce $p_j = |\mathcal{T}^+_j|$ and $n_j = |\mathcal{T}^-_j|$.

\section{Impact of Perfect Efficiencies, Positive Prices and Strictly Negative Prices} \label{sec:ideal}

We first formalize existing results for positive prices and perfect efficiencies. We further consider the possibility of strictly negative prices and introduce a refined MILP relaxation of the storage scheduling problem.

\subsection{Positive Prices and Perfect Efficiencies}

As seen in Section \ref{sec:related}, the combination of strictly positive prices and inefficiencies ($\eta <1$) has been identified as a sufficient condition for the suboptimality of simultaneous charge and discharge in various previous papers. 
In the case of perfect efficiency, relaxation \eqref{pb:opt} may have multiple optimal solutions, including some with simultaneous charge and discharge, but always at least one without. The same is valid for prices equal to zero, though the common assumption is that prices are strictly positive, which is unnecessarily restrictive. We formalize this in Proposition \ref{prop1}. 

\begin{proposition} \label{prop1}
    Some special cases for which relaxation \eqref{pb:opt} is exact are as follows:
    \begin{enumerate}
        \item If $\mathcal{T}^- \cup \mathcal{T}^0=\varnothing$ and $\eta < 1$, $\mathcal{X}\uptxt{LP*} \cap \mathcal{X}\uptxt{CC} = \mathcal{X}\uptxt{LP*}$.
        \item If $\eta = 1$, $\mathcal{X}\uptxt{LP*} \cap \mathcal{X}\uptxt{CC} \neq \varnothing$.
        \item If $\mathcal{T}^-=\varnothing$, $\mathcal{X}\uptxt{LP*} \cap \mathcal{X}\uptxt{CC} \neq \varnothing$.
    \end{enumerate}
\end{proposition}

Note that the proof of Proposition \ref{prop1}
contains the equations for obtaining an equivalent solution without simultaneous charge and discharge when $\eta = 1$ and when $C_t=0$.

We will now evaluate simultaneous charge and discharge when the strong conditions in Proposition~\ref{prop1} are not satisfied, which corresponds to Assumption \ref{ass}.
\begin{assumption}\label{ass}
   $\eta < 1$ and $\mathcal{T}^-\neq\varnothing$.
\end{assumption}

\subsection{Strictly Negative Prices and a Refined MILP Relaxation}

\begin{corollary} \label{corol1}
     Suppose that Assumption \ref{ass} is satisfied.
     If an optimal solution exhibits simultaneous charge and discharge in time period $t$, then $t \in \mathcal{T}^- \cup \mathcal{T}^0$. Moreover, if $t \in \mathcal{T}^0$, then there exists an alternative optimal solution without simultaneous charge and discharge in time period $t$.
\end{corollary} 

From Corollary \ref{corol1}, we can conclude that, in order to avoid simultaneous charge and discharge, it is sufficient to include complementarity constraints or their equivalent with binary variables in the hours with strictly negative prices, thereby reducing the computing burden. We therefore introduce a refined relaxation of \eqref{pb:opt_milp}, which is guaranteed to have an optimal solution that is also optimal for \eqref{pb:opt_milp}:
\begin{subequations}  \label{pb:opt_milp_2}
\begin{align}
    \label{eq:milp_2_obj} \max_{\mathbf{x}} \quad & Z \\
    \label{eq:milp_2_stg_mdl} \text{s.t.} \quad & \eqref{eq:update} - \eqref{eq:bounds_pc} \, , \\
    \label{eq:bounds_pc_nobin} & 0 \leq p_t\uptxt{C} \leq \overline{P}\uptxt{C} \, ,& \forall t \in \mathcal{T^+}, \\
    \label{eq:bounds_pd_nobin} & 0 \leq p_t\uptxt{D} \leq \overline{P}\uptxt{D} \, ,& \forall t \in \mathcal{T^+}, \\
    \label{eq:bounds_pc_bin_2} & 0 \leq p_t\uptxt{C} \leq u_t\uptxt{C} \overline{P}\uptxt{C} \, ,& \forall t \in \mathcal{T^-}, \\
    \label{eq:bounds_pd_bin_2} & 0 \leq p_t\uptxt{D} \leq u_t\uptxt{D} \overline{P}\uptxt{D} \, ,& \forall t \in \mathcal{T^-}, \\
    \label{eq:complementarity_bin_2} & u_t\uptxt{C} + u_t\uptxt{D} \leq 1 \, ,& \forall t \in \mathcal{T^-}, \\
    \label{eq:bin_def_2} & u_t\uptxt{C}, u_t\uptxt{D} \in \{ 0,1\} \, ,& \forall t \in \mathcal{T^-}.
\end{align}
\end{subequations}

However, even when there are strictly negative prices, simultaneous charge and discharge is not always optimal. This is the focus of Section \ref{sec:new}.

\subsection{Illustration}
To illustrate the results of Proposition \ref{prop1} and Corollary \ref{corol1}, we first consider the scheduling of a storage system with imperfect round-trip efficiency, and for positive prices. For this example and the following ones, time periods have a duration $\Delta t$ of one hour, there is no leakage ($\rho=1$), and unless stated otherwise, the horizon is 24 hours, the capacity bounds are $\underline{S}=0$ and $\overline{S}=1$ MWh, the initial state of energy is $S\uptxt{init}=0$, and efficiencies are $\eta\uptxt{C} = \eta\uptxt{D} = 0.9$. 
For this first example, we have $\overline{P}\uptxt{C} = \overline{P}\uptxt{D}=1.5$ MW.
\ifarxiv The code for all the examples presented in this paper is available online at \url{https://github.com/eleaprat/simultaneous_charge_discharge}. \fi

Here and in the rest of the paper, we schedule the storage system for the prices of the day-ahead market for the zone DK1 (Denmark) in 2023, a year with a high occurrence of strictly negative prices. In this first example, we use the prices of the $26\uptxt{th}$ of May 2023. They are shown in Figure \ref{fig:res_ex_1_1}, along with the results from scheduling the storage system. The prices on that day were all positive and one hour had a price of 0, which is indicated by the gray area on the plots. 

\ifarxiv
\begin{figure}[ht]
    \centering
    \subfloat[Result with simultaneous charge and discharge plotted against time (hours) \label{fig:res_ex_1_1_b}]
    {\resizebox{0.9\linewidth}{!}{\includegraphics{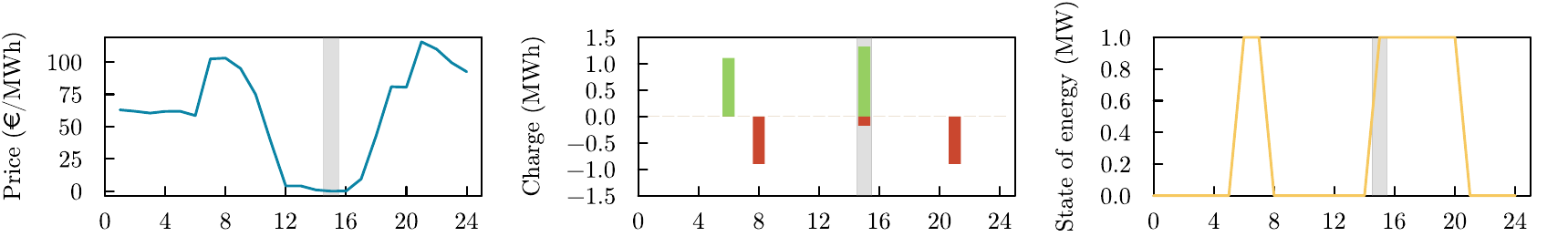}}}
    \\ \vspace{-0.2cm}
    \subfloat[Result without simultaneous charge and discharge plotted against time (hours) \label{fig:res_ex_1_1_a}]
    {\resizebox{0.9\linewidth}{!}{\includegraphics{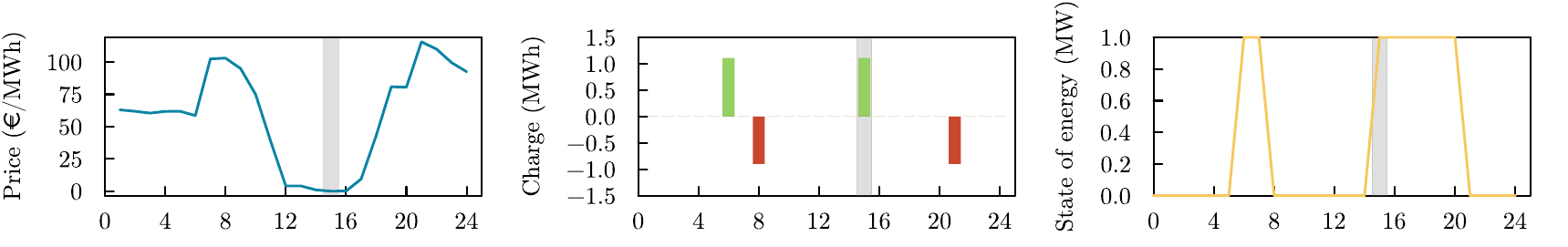}}}
    \caption{Two equivalent results for the illustration of Proposition \ref{prop1} and Corollary \ref{corol1}, with positive prices and $\eta <1$}
    \label{fig:res_ex_1_1}
\end{figure}

\begin{figure}[bht]
    \centering
    \subfloat[$\eta = 1$ \label{fig:res_ex_1_2}]
    {\resizebox{0.9\linewidth}{!}{\includegraphics{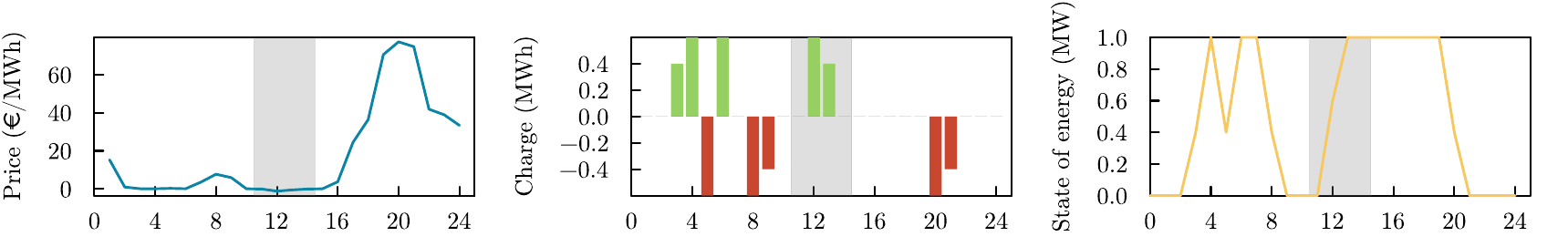}}}
    \\ \vspace{-0.4cm}
    \subfloat[$\eta < 1$ \label{fig:res_ex_1_3}]
    {\resizebox{0.9\linewidth}{!}{\includegraphics{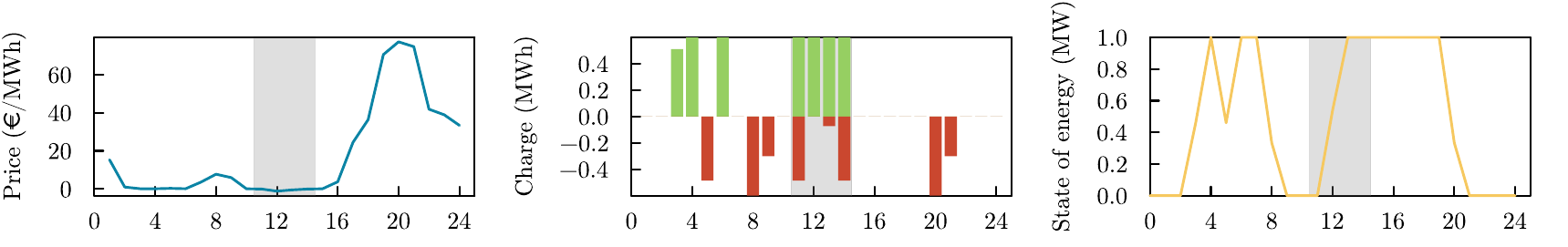}}}
    \\ \vspace{-0.4cm}
    \subfloat[$\eta < 1$, slower \label{fig:res_ex_1_4}]
    {\resizebox{0.9\linewidth}{!}{\includegraphics{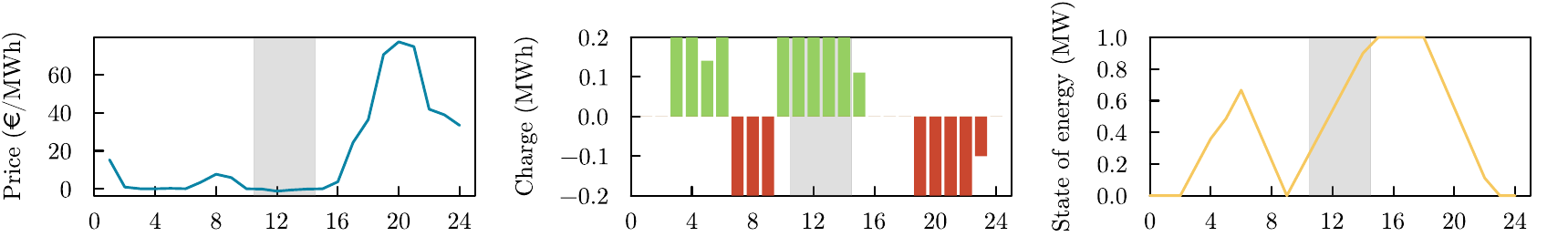}}}
    \caption{Results for the illustration of Proposition \ref{prop1} and Corollary \ref{corol1}, with strictly negative prices, plotted against time (hours)}
    \label{fig:res_ex_1_2_3}
\end{figure}
\fi

We observe that, as expected, there is no simultaneous charge and discharge when the price is strictly positive. When the price is zero, there can be simultaneous charge and discharge, which we see in Figure~\ref{fig:res_ex_1_1_b}; however, there is an equivalent solution without simultaneous charge and discharge, for which the state of energy is the same, in Figure \ref{fig:res_ex_1_1_a}.\footnote{Note that the solver returns the solution without simultaneous charge and discharge. To obtain the other solution, we add a constraint to impose discharge in that hour.} The profit is in both cases equal to 131.76 €.

\ifarxiv
\begin{table}[ht]
    \centering
    \caption{Data for the illustration of Proposition \ref{prop1} and Corollary \ref{corol1}, with strictly negative prices}
    \begin{tabular}{lrrr}
    \hline
    \textbf{Parameter}                         & \multicolumn{1}{c}{\textbf{Case $\eta = 1$}} & \multicolumn{1}{l}{\textbf{Case $\eta < 1$}} & \multicolumn{1}{l}{\textbf{Case $\eta < 1$, slower}} \\ \hline
    $\overline{P}\uptxt{C}$, $\overline{P}\uptxt{D}$ & 0.6 MW                                       & 0.6 MW                                       & 0.2 MW                                             \\
    $\eta\uptxt{C}$, $\eta\uptxt{D}$                 & 1                                            & 0.9                                          & 0.9                                                \\ \hline
    \end{tabular}
    \label{tab:data_ex_1_2}
\end{table}
\fi

To illustrate the results for strictly negative prices, we schedule the storage systems in Table \ref{tab:data_ex_1_2}, comparing results for $\eta = 1$ and for $\eta<1$.
We use the prices for the $25\uptxt{th}$ of March 2023. They are shown in Figure~\ref{fig:res_ex_1_2_3}, along with the results from scheduling the storage systems. On that day, four consecutive hours with strictly negative prices were observed, which is indicated by the gray area on the plots. 
We observe in Figure \ref{fig:res_ex_1_2} that there is no simultaneous charge and discharge when $\eta = 1$. Note that in this case, there also exist equivalent solutions with simultaneous charge and discharge.
When $\eta < 1$, we can see in Figure \ref{fig:res_ex_1_3} that there is simultaneous charge and discharge, and as we expected from Corollary \ref{corol1}, it is only in the hours of strictly negative prices. For one of these time periods, there is no simultaneous charge and discharge, which illustrates that simultaneous charge and discharge is not necessarily optimal, even for strictly negative prices.
We consider a last case for which we also have $\eta < 1$, but with a slower charge, which corresponds to the last column of Table \ref{tab:data_ex_1_2}. We can observe in Figure \ref{fig:res_ex_1_4}  that despite the strictly negative prices, there is no simultaneous charge and discharge. This example proves that there are cases when exclusivity of charge and discharge is ensured even in the presence of negative prices, showing that the conditions in Proposition~\ref{prop1} are not necessary.
In the following section, we aim to formalize the conditions under which this happens.

\section{New Conditions for Evaluating Simultaneous Charge and Discharge} \label{sec:new}

We first present an a posteriori result that will give us some insights regarding the link between simultaneous charge and discharge and the storage system characteristics. It is then used to prove our main results. We identify two sufficient conditions guaranteeing that any optimal solution of \eqref{pb:opt} will exhibit simultaneous charge and discharge.
We then give two other sufficient conditions that guarantee that the relaxation is exact.

\subsection{An a Posteriori Result}

Lemma \ref{lemma1} connects the optimality of simultaneous charge and discharge in problem \eqref{pb:opt} with the difference in the state of energy between two time periods.

\begin{lemma} \label{lemma1}
    For problem \eqref{pb:opt}, under Assumption \ref{ass}, considering $t \in \mathcal{T}^-$, if the optimal solution is to net charge the maximum quantity, i.e. $s_t^*-\rho s_{t-1}^* = \Delta t \, \eta\uptxt{C} \overline{P}\uptxt{C}$ (with $s_0^*=S\uptxt{init}$), or to net discharge the maximum quantity, i.e. $s_t^*-\rho s_{t-1}^* = - \Delta t ({1}/{\eta\uptxt{D}}) \overline{P}\uptxt{D}$, simultaneous charge and discharge is not optimal.
    In contrast, if this condition is not satisfied, simultaneous charge and discharge will occur at $t$ in any optimal solution.
\end{lemma} 

Lemma \ref{lemma1} demonstrates that having simultaneous charge and discharge is dependent on the ability to net charge or discharge the maximum quantity. The second part of this lemma indicates when an optimal solution will exhibit simultaneous charge and discharge. Indeed, if net charging or discharging the maximum quantity is not possible because of the characteristics of the storage, we can expect that solving \eqref{pb:opt} will return a physically infeasible solution. In the following, we use the results of Lemma \ref{lemma1} to relate the optimality of simultaneous charge and discharge to the characteristics of the storage system.

\subsection{First Sufficient Condition for Inexact Relaxation}

\subsubsection{Sufficient Condition.}

Lemma \ref{lemma1} shows that if net charging or discharging the maximum quantity is not possible because of the characteristics of the storage, the relaxation \eqref{pb:opt} would be inexact. In the next corollary, we give conditions based on the storage characteristics for which this is the case.

\begin{corollary} \label{corol2}
Supposing that Assumption \ref{ass} holds, the relaxation \eqref{pb:opt} will be inexact if the following conditions hold:
\begin{enumerate}
    \item $\rho \underline{S} + \Delta t \, \eta\uptxt{C} \overline{P}\uptxt{C} > \overline{S}$
    \item $\rho \overline{S} - \Delta t ({1}/{\eta\uptxt{D}}) \overline{P}\uptxt{D} < \underline{S}$
\end{enumerate}
\end{corollary} 

Note that without leakage ($\rho=1$), Condition 1 becomes $\Delta t \uptxt{C} < \Delta t$ and Condition 2 becomes $\Delta t \uptxt{D} < \Delta t$, which corresponds to saying that the duration of charge and the duration of discharge are both less than the duration of a time period in the model.

\subsubsection{Illustration.}

We consider an example with strictly negative prices and compare two cases. The data related to the two corresponding storage systems is shown in Table~\ref{tab:data_ex_2}. In the first case, the duration of charge is 0.56 hours and the duration of discharge is 0.45 hours, and the conditions of Corollary~\ref{corol2} hold.
In the second case, we consider a storage system with a slower charge and discharge, for which these conditions do not hold. 
We schedule the systems for the $1\uptxt{st}$ of January 2023, a day for which there was only one hour with a strictly negative price. 
Prices are plotted in Figure \ref{fig:res_ex_2}, along with the results. 

\ifarxiv
\begin{table}[ht]
    \centering
    \caption{Data for the illustration of Corollary \ref{corol2}}
    \begin{tabular}{lrr}
    \hline
    \textbf{Parameter}                         & \textbf{Fast storage}                     & \textbf{Slow storage}                     \\ \hline
    $\overline{P}\uptxt{C}$, $\overline{P}\uptxt{D}$ & 2 MW                                      & 0.4 MW                                    \\
    $\Delta t\uptxt{C}$                 & 0.56 h                                       & 2.77 h                                       \\
    $\Delta t\uptxt{D}$                 & 0.45 h                                   & 2.25 h                                      \\ \hline
    \end{tabular}
    \label{tab:data_ex_2}
\end{table}
\fi

We can see in Figure \ref{fig:res_ex_2_1} that for the storage system with fast charge and discharge, there is simultaneous charge and discharge, as expected from Corollary \ref{corol2}. We can draw a parallel with the results of Lemma~\ref{lemma1}: if the storage system were net charging the maximum, the bounds of the storage would be violated. Therefore, it can only net charge less than the maximum, and there is simultaneous charge and discharge in the hour with a strictly negative price.
In contrast, it is possible for the slow storage system to charge the maximum quantity without reaching the storage capacity limit, as shown in the right plot of Figure \ref{fig:res_ex_2_2}. As stated by Lemma \ref{lemma1}, there is no simultaneous charge and discharge in that hour for the slow storage system.

\ifarxiv
\begin{figure}[ht]
    \centering
    \subfloat[Fast storage \label{fig:res_ex_2_1}]
    {\resizebox{0.9\linewidth}{!}{\includegraphics{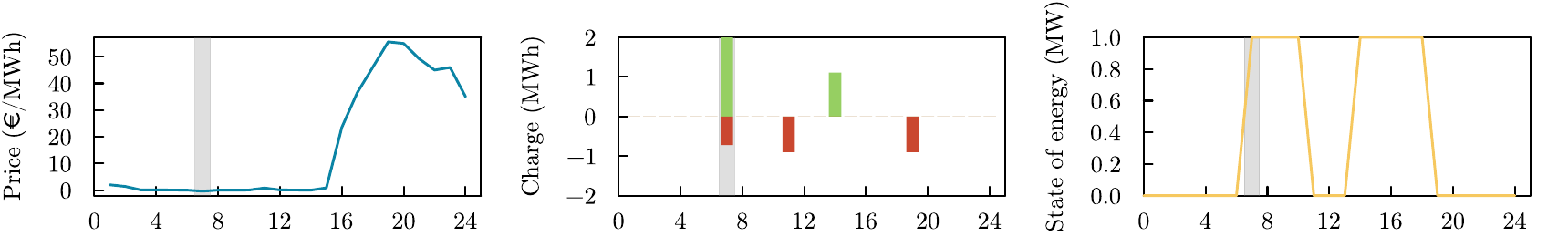}}}
    \\ \vspace{-0.4cm}
    \subfloat[Slow storage \label{fig:res_ex_2_2}]
    {\resizebox{0.9\linewidth}{!}{\includegraphics{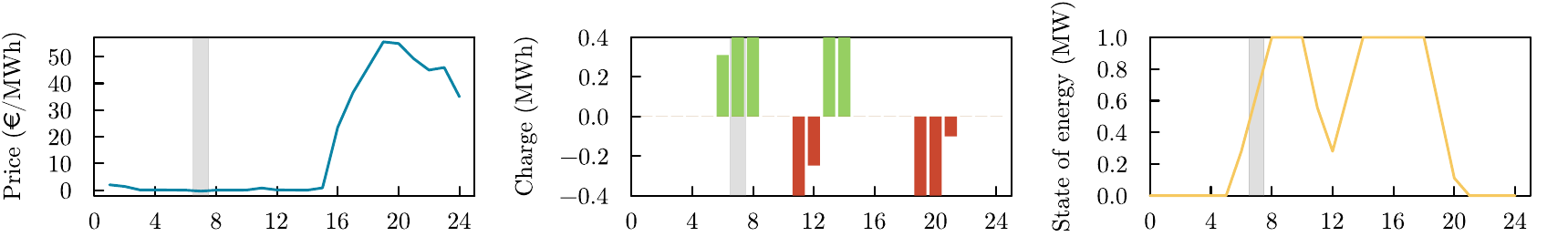}}}
    \caption{Results for the illustration of Corollary \ref{corol2}, plotted against time (hours)}
    \label{fig:res_ex_2}
\end{figure}
\fi

\subsection{Second Sufficient Condition for Inexact Relaxation}

\subsubsection{Sufficient Condition.}

We give a more general condition on the storage characteristics for which any optimal solution of \eqref{pb:opt} will exhibit simultaneous charge and discharge in the following theorem, which compares the time it takes for the storage system to fully charge to the duration of the longest sequence of negative prices, $\overline{n}$. 

\begin{theorem} \label{theorem1}
Suppose that Assumptions \ref{ass3} and \ref{ass} hold. The relaxation \eqref{pb:opt} will be inexact if the following conditions hold:
\begin{enumerate}
    \item $\rho^{\overline{n}} \underline{S} + \Delta t  \eta \uptxt{C} \overline{P}\uptxt{C} (\sum_{t=1}^{\overline{n}} \rho^{\overline{n}-t}) > \overline{S} $
    \item $\{\mathcal{T}\uptxt{C} \subseteq \mathcal{T}, \, \mathcal{T}\uptxt{D} \subseteq \mathcal{T}, \, s\in \mathbb{R}^+ | \, \mathcal{T} \uptxt{C} \cup \mathcal{T} \uptxt{D} = \overline{\mathcal{T}}^-, \, \underline{S} \leq s \leq \overline{S}, $ \\ $  \qquad \qquad  \qquad \, \rho^{\overline{n}} s + \Delta t ( \eta \uptxt{C} \overline{P}\uptxt{C} ( \sum_{t \in \mathcal{T} \uptxt{C}} \rho^{\tau_2-t}) - ({1}/{\eta \uptxt{D}}) \overline{P}\uptxt{D} ( \sum_{t \in \mathcal{T} \uptxt{D}} \rho^{\tau_2-t}) ) = \overline{S} \} = \varnothing$, where $\overline{\mathcal{T}}^- = [\tau_1, \tau_2]$
\end{enumerate}
\end{theorem} 

With Theorem \ref{theorem1}, we can quickly identify cases where simultaneous charge and discharge can be a problem. Condition 1 corresponds to saying that the storage system can fully charge in fewer time periods than~$\overline{n}$. Condition 2 is harder to interpret, so we illustrate it with an example in what follows. Though this condition may seem restrictive, it is actually only relevant for the actual storage system level at the beginning of this longest sequence of strictly negative prices, i.e., for $s=s_{\tau_1-1}^*$.
Fixing the value of $s$ to $s_{\tau_1-1}^*$ considerably reduces the number of solutions of the system in Condition 2. 
Therefore, for a given problem, if Condition 1 holds, we recommend solving \eqref{pb:opt_milp_2} instead of \eqref{pb:opt}. 

\subsubsection{Illustration.}

To illustrate the results of Theorem \ref{theorem1}, we identify the longest sequence of strictly negative prices in 2023 for DK1. It consists of 32 hours of negative prices between the $7\uptxt{th}$ and the $9\uptxt{th}$ of July 2023, shown in the first plot of Figure \ref{fig:res_ex_3_1} and Figure \ref{fig:res_ex_3_2}. Therefore, $\overline{n}=32$. We consider three different storage systems, with the characteristics presented in Table~\ref{tab:data_ex_3}. 
Without leakage, Condition 1 in Theorem \ref{theorem1} reduces to $ \underline{S} + \Delta t  \eta \uptxt{C} \overline{P}\uptxt{C} \overline{n} > \overline{S} $. For our capacity bounds, the critical value of the maximum charge is thus $\overline{P}\uptxt{C} = 0.035$ MW. 

\ifarxiv
\begin{table}[ht]
    \centering
    \caption{Data for the illustration of Theorem \ref{theorem1}}
    \begin{tabular}{lrrr}
    \hline
    \textbf{Parameter}         & \textbf{Case 1}                     & \textbf{Case 2}                     & \textbf{Case 3}                     \\ \hline
    $\overline{P}\uptxt{C}$       & 0.036 MW                            & 0.034 MW                            & 0.036 MW                            \\
    $\overline{P}\uptxt{D}$       & 0.036 MW                            & 0.028 MW                            & 0.00396 MW                          \\
    $\Delta t \uptxt{C}$ & 30.86 h                                 & 32.68 h                                 & 30.86 h                                 \\
    $\Delta t \uptxt{D}$ & 25 h                                 & 32.14 h                                 & 227.27 h                                \\ \hline
    \end{tabular}
    \label{tab:data_ex_3}
\end{table}
\fi

In the first case, we consider $\overline{P}\uptxt{C} > 0.035$ MW, so Condition 1 holds.
We can see in Figure \ref{fig:res_ex_3_1} that there is simultaneous charge and discharge in the first and last hours of this sequence of strictly negative prices. The storage system is scheduled from the initial level $S\uptxt{init}=0$, which corresponds to $s_{\tau_1-1}^*$ here. We do verify that the equations in Condition 2 do not have a solution for $s=0$. Therefore, the solution observed, with simultaneous charge and discharge, is as expected from Theorem \ref{theorem1}. The storage system reaches its capacity in less than 32 hours; therefore, it is unable to net charge the maximum quantity in all the hours of the longest sequence of negative prices. Consequently, simultaneous charge and discharge is optimal.

\ifarxiv
\begin{figure}[ht]
    \centering
    \subfloat[Case 1 \label{fig:res_ex_3_1}]
    {\resizebox{0.9\linewidth}{!}{\includegraphics{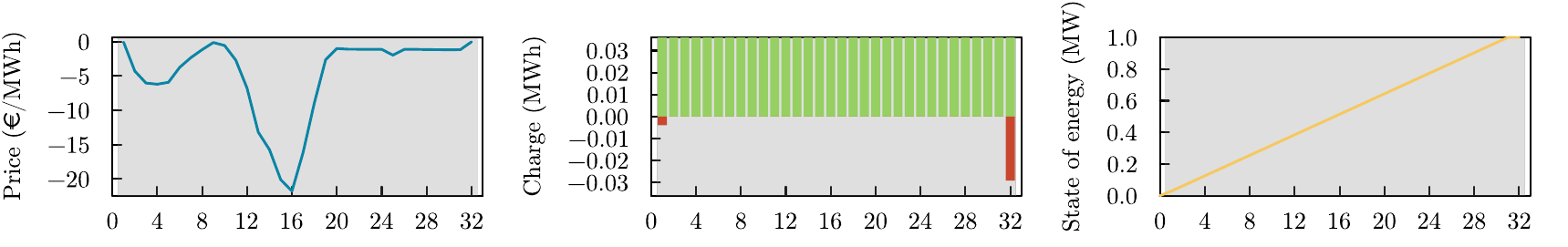}}}
    \\ \vspace{-0.4cm}
    \subfloat[Case 2 \label{fig:res_ex_3_2}]
    {\resizebox{0.9\linewidth}{!}{\includegraphics{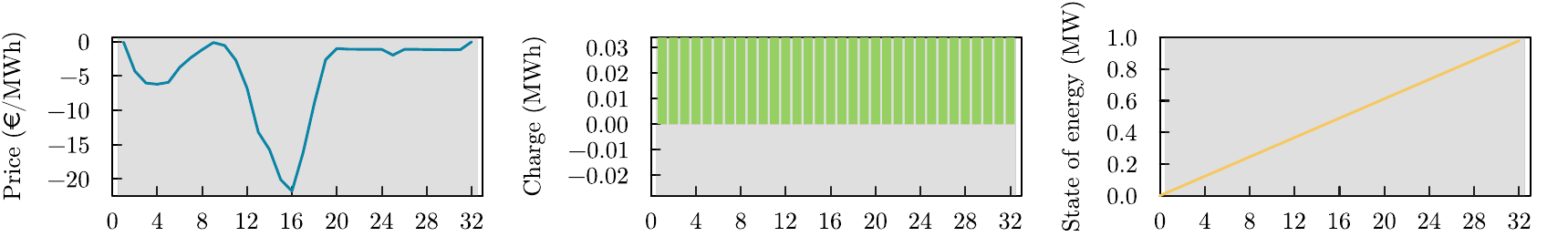}}}
    \\ \vspace{-0.4cm}
    \subfloat[Case 3 \label{fig:res_ex_3_3}]
    {\resizebox{0.9\linewidth}{!}{\includegraphics{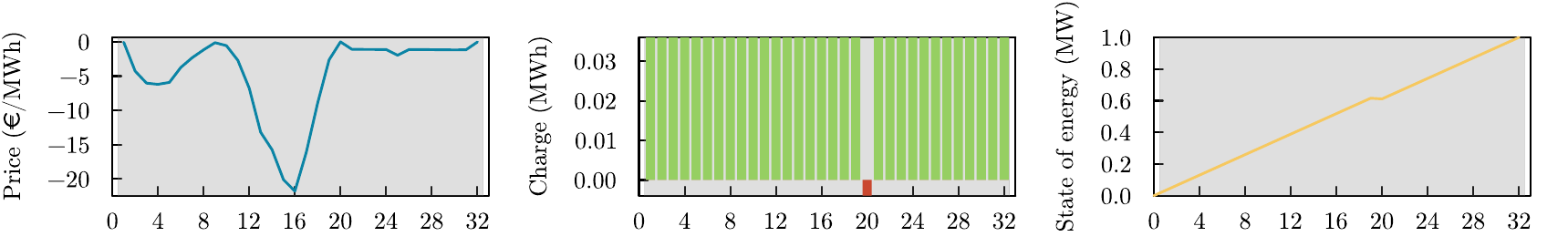}}}
    \caption{Results for the illustration of Theorem \ref{theorem1}, plotted against time (hours), considering the longest sequence of strictly negative prices}
    \label{fig:res_ex_3}
\end{figure}
\fi

In contrast, we can see the results from scheduling a storage system with $\overline{P}\uptxt{C} \leq 0.035$ MW in Figure~\ref{fig:res_ex_3_2}. In this case, even though there may not be simultaneous charge and discharge when scheduling for each sequence of strictly negative prices individually, there is no guarantee that this will still be the case when they are considered in the same horizon. This corresponds to the situation illustrated in Figure \ref{fig:res_ex_3_long_2}, which includes the 11 previous hours, in which there is another sequence of strictly negative prices. In this case, the duration between the two sequences is too short for the storage system to discharge enough so that it can net charge the maximum quantity in all the hours with strictly negative prices.\footnote{Note that the speed of discharge is slower than the speed of charge.} As a consequence, there is simultaneous charge and discharge in the first hour. This shows that Theorem \ref{theorem1} is useful for identifying storage system characteristics for which simultaneous charge and discharge will occur, in the example, for $\overline{P}\uptxt{C} > 0.035$~MW, but it does not cover all the cases of simultaneous charge and discharge.

\ifarxiv
\begin{figure}[ht]
    \centering
    \subfloat[Case 2 \label{fig:res_ex_3_long_2}]
    {\resizebox{0.9\linewidth}{!}{\includegraphics{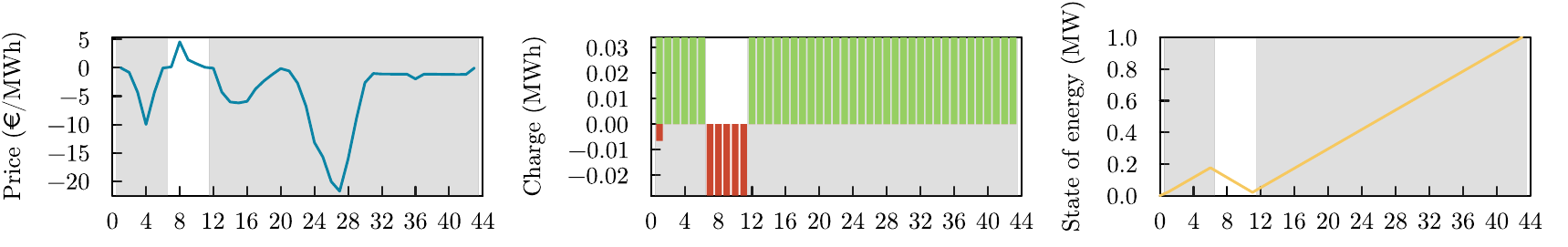}}}
    \\ \vspace{-0.4cm}
    \subfloat[Case 3 \label{fig:res_ex_3_long_3}]
    {\resizebox{0.9\linewidth}{!}{\includegraphics{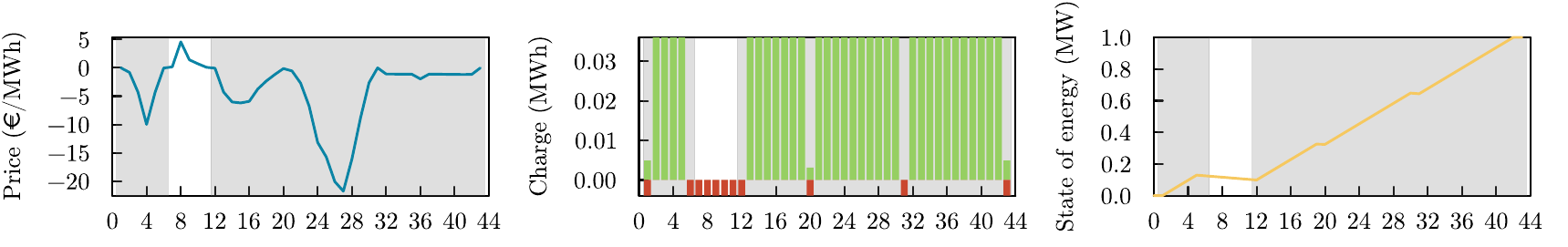}}}
    \caption{Results for the illustration of Theorem \ref{theorem1}, plotted against time (hours), considering two close sequences of strictly negative prices}
    \label{fig:res_ex_3_long}
\end{figure}
\fi

Condition 2 of Theorem \ref{theorem1} is necessary to cover the unlikely case in which the storage system would violate the capacity bound if it were to net charge the maximum quantity in all the time periods of this longest sequence of strictly negative prices, but could \emph{exactly} reach the bound by net discharging the maximum quantity in some of these time periods. In this case, the sequence of strictly negative prices is so long compared to the storage size that the storage system is able to arbitrage between the hours of negative prices. This can only happen if the time periods with the highest prices are in the middle of the sequence of negative prices, which is also unlikely in practice. For the time series of $7\uptxt{th}$ to $9\uptxt{th}$ July 2023, the highest prices are in the first and in the last hours, which is why there is simultaneous charge and discharge in these hours in the first case, shown on Figure \ref{fig:res_ex_3_1}. In the first time period, net discharging is not possible because it would violate the lower bound on the storage capacity. Moreover, it is necessary to net discharge before the last hour, otherwise the upper bound on capacity would be violated. Therefore, to illustrate Condition~2, we modify the sequence of prices so that the highest price is in the middle, see Figure \ref{fig:res_ex_3_3} (in hour 20, the price is modified to -0.01 €/MWh). We schedule the storage system with the characteristics given in Table~\ref{tab:data_ex_3} for Case 3. For this storage system, the equations of Condition 2 have a solution for $s = S\uptxt{init} = 0$, with $|\mathcal{T}\uptxt{D}| = 31$ and $|\mathcal{T}\uptxt{C}|=1$. This means that it is possible for this storage system, starting from the lowest level, to reach the maximum level by net charging the maximum quantity during 31 hours and net discharging the maximum quantity during 1 hour. This solution is preferable compared to net charging or net discharging less than the maximum quantities, which can be observed in Figure \ref{fig:res_ex_3_3}. In this case, there is no simultaneous charge and discharge for this sequence of strictly negative prices. However, this corresponds to very specific characteristics and conditions.
In fact, if we include the 11 previous hours, we can see in Figure \ref{fig:res_ex_3_long_3} that there is simultaneous charge and discharge.
This demonstrates that even when Conditions 1 and 2 stand, there is no guarantee that simultaneous charge and discharge will not occur, and solving \eqref{pb:opt_milp_2} should be preferred.

\subsection{First Sufficient Condition for Exact Relaxation}

\subsubsection{Sufficient Condition.}

We now focus on the case where there is only one sequence of strictly negative prices. In this setup, we give sufficient conditions under which the relaxation is exact.

\begin{theorem} \label{theorem2}
    Suppose that Assumption \ref{ass} holds. Under the following conditions, relaxation \eqref{pb:opt} is exact:
    \begin{enumerate}
        \item $\mathcal{T} = \mathcal{T}^+_1 \cup \mathcal{T}^-_1 \cup \mathcal{T}^+_2 $, with $\mathcal{T}^-_1 = [\tau_1,\tau_2] \neq \varnothing$.
        \item $ \rho^{\tau_2-\tau_1} ( \max \{\rho^{\tau_1-1} S\uptxt{init}- \Delta t ({1}/{\eta\uptxt{D}})\overline{P}\uptxt{D} (\sum_{t=1}^{\tau_1-1}{\rho^{\tau_1-1-t}}) ,\underline{S} \} ) + \Delta t \, \eta\uptxt{C} \overline{P}\uptxt{C} (\sum_{t=\tau_1}^{\tau_2}\rho^{\tau_2-t}) \leq \overline{S} $
    \end{enumerate}
\end{theorem}

The idea behind Theorem \ref{theorem2} is that if the storage system can charge at the maximum rate, 
there will not be simultaneous charge and discharge in the time periods with a strictly negative price. The second condition essentially states that full charge is possible, i.e, that the capacity of the storage is not exceeded when net charging the maximum quantity between $\tau_1$ and $\tau_2$, starting from the initial level and considering maximum discharge until $\tau_1$, or until reaching the minimum level. In particular, for a storage system without leakage, this is the case if it is initially at its minimum level and its duration of charge satisfies $\Delta t \uptxt{C} \geq \tau_2-\tau_1+1$.

\subsubsection{Illustration.}

To illustrate the results of Theorem \ref{theorem2}, we consider that there is only one occurrence of strictly negative prices, and for this, we use again the prices of $25\uptxt{th}$ March 2023. There is exactly one interval of time periods for which the prices are strictly negative, with $\tau_1 = 11$ and $\tau_2 = 14$. The first 10 hours of positive prices are followed by 4 hours of strictly negative prices. 
For a storage system without leakage, the last condition thus becomes $ \max \{S\uptxt{init} - 10 \Delta t ({1}/{\eta\uptxt{D}})\overline{P}\uptxt{D},\underline{S} \} + 4  ( \Delta t \, \eta\uptxt{C} \overline{P}\uptxt{C}  )\leq \overline{S} $: starting from the initial level, if discharging as much as possible during 10 hours, and charging as much as possible during 4 hours, the maximum level is not exceeded.
We check how this condition is impacted by different storage characteristics and initial levels. 
The different cases are summarized in Table \ref{tab:data_ex_4}.

\ifarxiv
\begin{table}[ht]
    \centering
    \caption{Data for the illustration of Theorem \ref{theorem2}}
    \begin{tabular}{lrrr}
    \hline
    \textbf{Parameter}         & \textbf{Case 1}                     & \textbf{Case 2}                     & \textbf{Case 3}                     \\ \hline
    $\overline{P}\uptxt{C}$       & 0.27 MW                             & 0.27 MW                             & 0.28 MW                            \\
    $\overline{P}\uptxt{D}$       & 0.27 MW                             & 0.08 MW                             & 0.27 MW                          \\
    $\Delta t \uptxt{C}$ & 4.12 h                                 & 4.12 h                                 & 3.97 h                                 \\
    $\Delta t \uptxt{D}$ & 3.33 h                                 & 11.25 h                                 & 3.33 h                                \\ \hline
    \end{tabular}
    \label{tab:data_ex_4}
\end{table}
\fi

First, we consider that $S\uptxt{init} = \underline{S}$ and $\overline{P}\uptxt{C} \leq ({\overline{S}-\underline{S}})/({4 \Delta t \, \eta\uptxt{C}})$, where $({\overline{S}-\underline{S}})/({4 \Delta t \, \eta\uptxt{C}}) = 0.278$ MW, so that Condition 2 holds. The results are shown in Figure \ref{fig:res_ex_4_1}. Note that Condition 2 does not correspond to the actual optimal scheduling, but rather checks that i) it is possible to discharge enough before the sequence of negative prices, and ii) charge is not too fast. Here we can see that the storage actually performs some arbitrage before the sequence of strictly negative prices and even charges a small quantity, such that the maximum level is exactly attained at the end of this sequence.

\ifarxiv
\begin{figure}[ht]
    \centering
    \subfloat[Case 1 \label{fig:res_ex_4_1}]
    {\resizebox{0.9\linewidth}{!}{\includegraphics{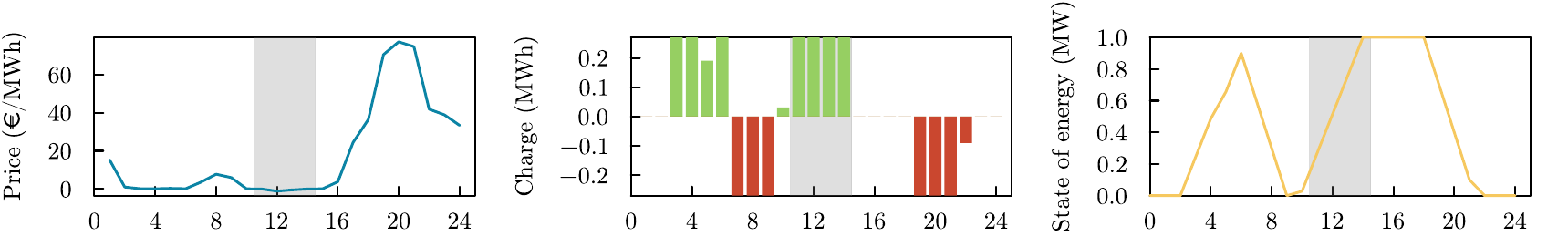}}}
    \\ \vspace{-0.4cm}
    \subfloat[Case 2 \label{fig:res_ex_4_2}]
    {\resizebox{0.9\linewidth}{!}{\includegraphics{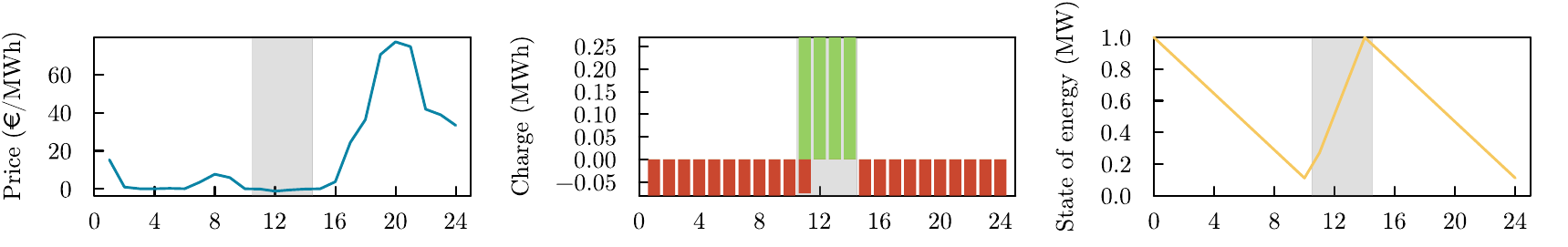}}}
    \\ \vspace{-0.4cm}
    \subfloat[Case 3 \label{fig:res_ex_4_3}]
    {\resizebox{0.9\linewidth}{!}{\includegraphics{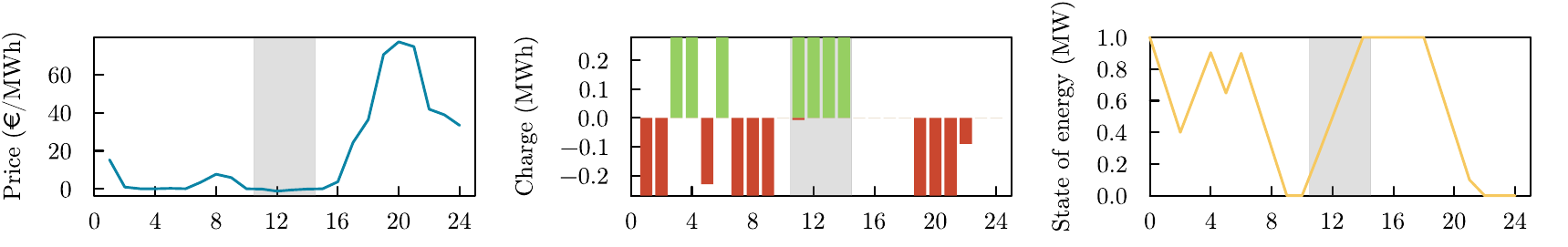}}}
    \caption{Results for the illustration of Theorem \ref{theorem2}, plotted against time (hours)}
    \label{fig:res_ex_4}
\end{figure}
\fi

In the second case, we increase the initial level and decrease the discharging rate, so that the minimum level cannot be reached before the sequence of strictly negative prices. The results are in Figure \ref{fig:res_ex_4_2}. In this case, though the storage system net discharges the maximum quantity until the sequence of strictly negative prices, net charging is not possible in all of the hours of this sequence. As a result, simultaneous charge and discharge is optimal. 
This illustrates the influence of the initial level and of the discharging rate.

In the third case, the initial level and the discharging rate of Case 1 are kept but the charging rate is increased. The results are in Figure \ref{fig:res_ex_4_3}. The lowest level can be reached as in the first case, but net charging the maximum quantity during four hours is not possible. Note that the charging speed is only slightly higher than the critical value defined by Condition 2, and the quantity discharged in the hour of simultaneous charge and discharge is very small. This suggests that the shorter the duration of charge, the further we can find ourselves from the optimal solution of \eqref{pb:opt_nl} when solving \eqref{pb:opt}. In other words, if the storage system duration of charge is long enough, the effects of simultaneous charge and discharge might be negligible.

\subsection{Second Sufficient Condition for Exact Relaxation}

\subsubsection{Sufficient Condition.}

The next theorem gives the most general condition, which can be applied if we are not able to conclude from any of the previous results.

\begin{theorem} \label{theorem3}
Suppose that Assumptions \ref{ass3} and \ref{ass} hold. Considering that $\mathcal{T} = \bigcup_{j \in \mathcal{J}} ( \mathcal{T}^+_j \cup \mathcal{T}^-_j )$, relaxation~\eqref{pb:opt} is exact under the following conditions:
\begin{enumerate}
    \item $\hat{S}_j \leq \overline{S}$, $\forall j \in \mathcal{J}$
    \item $\hat{S}_j = \rho^{n_j} \max \{ \rho^{p_j} \hat{S}_{j-1} - \Delta t ({1}/{\eta\uptxt{D}}) \overline{P}\uptxt{D} (\sum_{t=1}^{p_j} \rho^{p_j-t}), \underline{S}\} + \Delta t \, \eta\uptxt{C} \overline{P}\uptxt{C} (\sum_{t=1}^{n_j} \rho^{n_j-t})$, $\forall j \in \mathcal{J} \setminus\{1\}$
    \item $\hat{S}_1 = \rho^{n_1} \max \{ \rho^{p_1} S\uptxt{init} - \Delta t ({1}/{\eta\uptxt{D}}) \overline{P}\uptxt{D} (\sum_{t=1}^{p_1} \rho^{p_1-t}), \underline{S}\} + \Delta t \, \eta\uptxt{C} \overline{P}\uptxt{C}(\sum_{t=1}^{n_1} \rho^{n_1-t})$.
\end{enumerate}
\end{theorem}

Theorem \ref{theorem3} says that there is a solution without simultaneous charge and discharge if, for each sequence of strictly negative prices, starting from the lowest reachable level, when net charging the maximum quantity, the maximum level is not exceeded. It implies that the storage system can discharge sufficiently between two sequences of strictly negative prices.

\subsubsection{Illustration.}

In the following, we show how to check the conditions of Theorem \ref{theorem3} on two examples: one with a storage system with fast charge and discharge, and one with a storage system with slow charge and discharge, as shown in Table \ref{tab:data_ex_5}. We use the prices from $27\uptxt{th}$ to $29\uptxt{th}$ May 2023, with three sequences of strictly negative prices, see Figure \ref{fig:res_ex_5}.

\ifarxiv
\begin{table}[hb]
    \centering
    \caption{Data for the illustration of Theorem \ref{theorem3}}
    \begin{tabular}{lrr}
    \hline
    \textbf{Parameter}                         & \textbf{Fast storage}                     & \textbf{Slow storage}                     \\ \hline
    $\overline{P}\uptxt{C}$, $\overline{P}\uptxt{D}$ & 0.2 MW                                    & 0.1 MW                                    \\
    $\Delta t \uptxt{C}$                 & 5.56 h                                     & 11.11 h                                       \\ 
    $\Delta t \uptxt{D}$                 & 4.5 h                                       & 9 h                                      \\ \hline
    \end{tabular}
    \label{tab:data_ex_5}
\end{table}
\fi

The successive values of $\hat{S}_j$ are shown in Table \ref{tab:res_ex_5}. In the case of the fast storage system, we can see that $\hat{S}_2 > \overline{S}$, meaning that the storage system cannot discharge sufficiently before the second sequence of strictly negative prices to be able to net charge the maximum quantity in all the hours of this sequence of strictly negative prices. Condition~1 is not satisfied and there is no need to continue the calculation for the next values of $j$. We can see in Figure \ref{fig:res_ex_5_1} that there is simultaneous charge and discharge in this case. For the slow storage system, on the other hand, all the values satisfy $\hat{S}_j \leq \overline{S}$ and all the conditions of Theorem~\ref{theorem3} are satisfied. We verify in Figure \ref{fig:res_ex_5_2} that there is no simultaneous charge and discharge in this case.

\ifarxiv
\begin{table}[ht]
    \centering
    \caption{Results for the illustration of Theorem \ref{theorem3}}
    \begin{tabular}{rrrrr}
    \hline
    \multicolumn{1}{l}{\textbf{$j$}} & \multicolumn{1}{l}{\textbf{$p_j$}} & \multicolumn{1}{l}{\textbf{$n_j$}} & \multicolumn{1}{l}{\textbf{$\hat{S}_j$ - Fast storage}} & \multicolumn{1}{l}{\textbf{$\hat{S}_j$ - Slow storage}} \\ \hline
    1                                & 13                                 & 2                                   & 0.36                                                    & 0.18\\
    2                                & 19                                 & 8                                   & \textbf{1.44}                                           & 0.72\\
    3                                & 16                                 & 8                                   & -                                                       & 0.72\\
    4                                & 6                                  & 0                                   & -                                                       & 0.05 \\ \hline
    \end{tabular}
    \label{tab:res_ex_5}
\end{table}
\fi

\ifarxiv
\begin{figure}[thb]
    \centering
    \subfloat[Fast storage \label{fig:res_ex_5_1}]
    {\resizebox{1.0\linewidth}{!}{\includegraphics{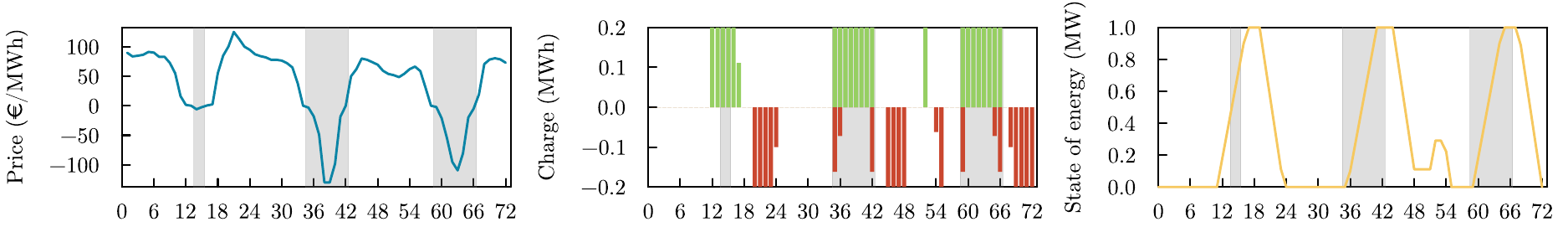}}}
    \\ \vspace{-0.4cm}
    \subfloat[Slow storage \label{fig:res_ex_5_2}]
    {\resizebox{1.0\linewidth}{!}{\includegraphics{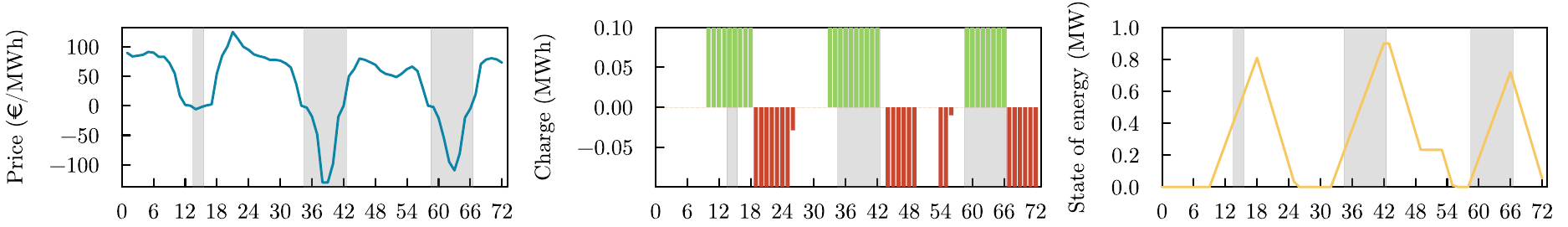}}}
    \caption{Results for the illustration of Theorem \ref{theorem3}, plotted against time (hours)}
    \label{fig:res_ex_5}
\end{figure}
\fi

\section{Discussion} \label{sec:discussion}

In this section, we show how the different results can be used to know if the linear relaxation of the storage system scheduling problem \eqref{pb:opt} will return the optimal solution of the original problem \eqref{pb:opt_nl}, or if the MILP formulation \eqref{pb:opt_milp_2} should be used. We also discuss some limitations and possible extensions.

The flowchart in Figure \ref{fig:dt} shows how to proceed when presented with a storage scheduling problem.
First, according to Proposition \ref{prop1}, in the case of perfect efficiency or positive prices only, \eqref{pb:opt} can be solved.
We then have the result of Corollary \ref{corol2}: if the storage system can charge completely and discharge completely in less than one time period, binary variables are needed to avoid simultaneous charge and discharge.
If not, we evaluate if there is a constraint on the final level of the storage system. In the models proposed here, in \eqref{pb:opt_nl}, its MILP reformulation \eqref{pb:opt_milp_2}, and in the linear relaxation \eqref{pb:opt}, we did not include the possibility of having a constraint imposing the final level of the storage system. This would actually add complexity to the results, as we would also need to consider that the final level should be reachable. It can increase the occurrence of simultaneous charge and discharge when solving \eqref{pb:opt}, and for this reason, we suggest solving \eqref{pb:opt_milp_2} if it is the case, also see \cite{Arroyo2022Ensuring}. Terminal conditions could also be included in the form of a penalty in the objective function. The impact of imposing a final level should be studied in future work. Note that when solving the scheduling problem using a rolling-horizon approach, in which the decisions of the first periods only are implemented while the rest of the horizon is advisory, the role of guiding the final level is ensured by including a longer horizon, such that an end-of-horizon level is not necessary.
If the final level is free, we can identify the longest sequence of strictly negative prices. If the storage system can charge completely in less than this duration, \eqref{pb:opt_milp_2} should be solved.
Otherwise, the results of either Theorem \ref{theorem2} or of Theorem \ref{theorem3} can be used, depending on whether there is one or more sequences of strictly negative prices.

\ifarxiv
\begin{figure}[htb]
    \centering
    {\resizebox{.7\linewidth}{!}{\includegraphics{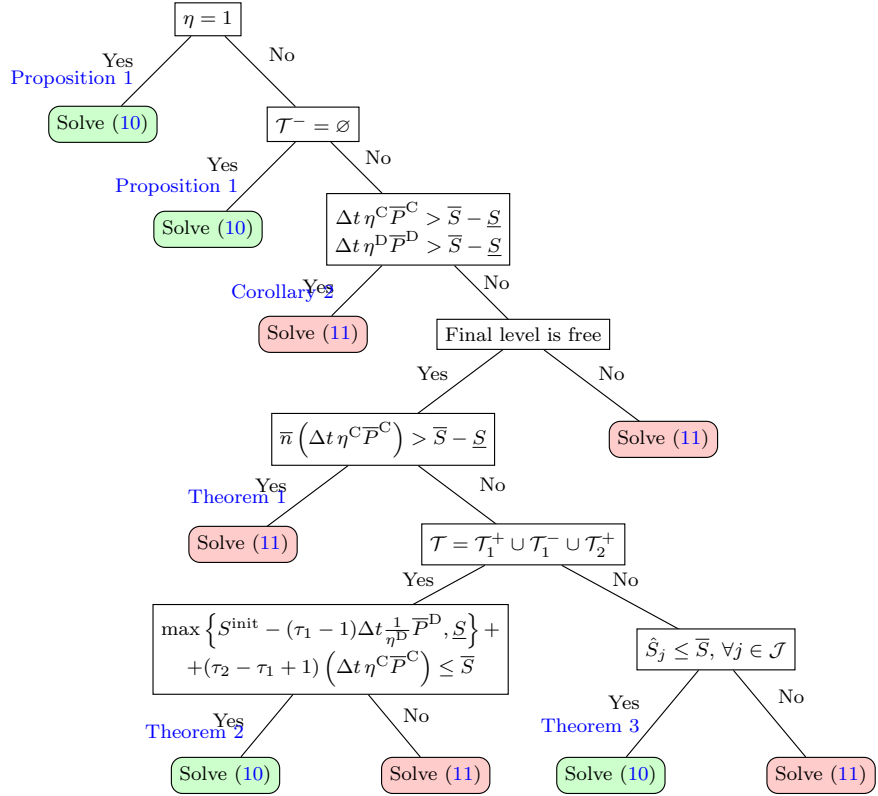}}}
    \caption{Flowchart describing how to use the results of this paper}
    \label{fig:dt}
\end{figure}
\fi

\section{Conclusion} \label{sec:ccl}
In this paper, we review and extend the list of a priori conditions for the exactness of the relaxation of the complementarity constraints in the energy storage scheduling problem. 
Conditions that were identified previously to this work are the perfect effectiveness of the storage system and the positivity of the prices received. As these assumptions do not correspond to reality, we propose new conditions that are less restrictive. We rigorously define these conditions, prove them mathematically, and illustrate them with examples. They can be summarized under the following assumption: the storage system duration of charge is long enough compared to the occurrence of negative prices. 
If this is not the case, it might be preferable to solve the MILP equivalent of the problem. We also introduce a refined formulation for this MILP equivalent, in which the number of binary variables is significantly reduced.
A natural extension would be to study the relevance of these results when considering the scheduling of storage systems as part of other problems, such as economic dispatch or optimal power flow, or for the storage of other commodities.
Indeed, while this study is motivated by the challenges of energy storage, the model used and the conclusions are valid for any type of storage.
Finally, in one of the examples, we observe that if the storage system duration of charge is slightly shorter than the critical value for which there is no simultaneous charge and discharge, the quantity affected might be negligible.
Therefore, future work could focus on analyzing the trade-off between computational complexity and the loss of optimality from solving the linear relaxation when the conditions identified in this paper are not satisfied, depending on the characteristics of the storage system. 

\section*{Code and Data}
The code and data for all the examples presented in this paper are available online at \url{https://github.com/eleaprat/simultaneous_charge_discharge}. 

\begin{appendices}
    
\section{Proof Proposition \ref{prop1}} \label{app:prop1}
First, we introduce the dual variables for \eqref{pb:opt}, repeated here for convenience: 
\begin{subequations} \label{pb:opt_dual}
\allowdisplaybreaks
\begin{align} 
    \label{eq:dual_obj} \max_{p\uptxt{D}, p\uptxt{C}, s} \quad & \sum_{t \in \mathcal{T}} \Delta t \, C_t (p_t\uptxt{D}-p_t\uptxt{C}) \\
    \label{eq:dual_init} \text{s.t.} \quad & s_1 = \rho S\uptxt{init} +  \Delta t (\eta\uptxt{C} p_1\uptxt{C} - \frac{1}{\eta\uptxt{D}} p_1\uptxt{D}) \, : \lambda_1 \\
    \nonumber & s_t = \rho s_{t-1} +  \Delta t 
     (\eta\uptxt{C} p_t\uptxt{C} - \frac{1}{\eta\uptxt{D}} p_t\uptxt{D}) \, , \\
     \label{eq:dual_update} & \qquad \qquad \qquad \qquad \quad \forall t \in \mathcal{T} \setminus \{1\} \, :\lambda_t\\
    \label{eq:dual_bounds_s} & \underline{S} \leq s_t \leq \overline{S} \, , \forall t \in \mathcal{T} \, :\underline{\sigma}_t,\overline{\sigma}_t \\
    \label{eq:dual_bounds_pc} & 0 \leq p_t\uptxt{C} \leq \overline{P}\uptxt{C} \, , \forall t \in \mathcal{T} \, :\underline{\gamma}_t,\overline{\gamma}_t\\
    \label{eq:dual_bounds_pd} & 0 \leq p_t\uptxt{D} \leq \overline{P}\uptxt{D} \, , \forall t \in \mathcal{T} \, :\underline{\delta}_t,\overline{\delta}_t.
\end{align}
\end{subequations}

The dual variable for each constraint is given after the colon on the right of the constraint.

The Karush-Kuhn--Tucker optimality conditions (KKTs) for \eqref{pb:opt} are:
\begin{subequations} \label{pb:kkts}
    \allowdisplaybreaks
    \begin{align}
        \label{eq:kkt1} & - \Delta t \, C_t + \lambda_t \frac{\Delta t}{\eta\uptxt{D}} + \overline{\delta}_t - \underline{\delta}_t = 0 \, ,&& \forall t \in \mathcal{T} \\
        \label{eq:kkt2} & \Delta t \, C_t - \lambda_t \Delta t \, \eta\uptxt{C} + \overline{\gamma}_t - \underline{\gamma}_t = 0  \, ,&& \forall t \in \mathcal{T} \\
        \label{eq:kkt3} & \lambda_t - \rho \lambda_{t+1} + \overline{\sigma}_t - \underline{\sigma}_t = 0  \, ,&& \forall t \in \mathcal{T} \\
        \label{eq:kkt4} & s_1 = \rho S\uptxt{init} +  \Delta t (\eta\uptxt{C} p_1\uptxt{C} - \frac{1}{\eta\uptxt{D}} p_1\uptxt{D})\\
        \label{eq:kkt5} & s_t = \rho s_{t-1} +  \Delta t (\eta\uptxt{C} p_t\uptxt{C} - \frac{1}{\eta\uptxt{D}} p_t\uptxt{D}) \, ,&& \, \forall t \in \mathcal{T} \setminus \{1\}\\
        \label{eq:kkt6} & 0 \leq s_t - \underline{S} \perp \underline{\sigma}_t \geq 0 \, ,&& \forall t \in \mathcal{T}\\
        \label{eq:kkt7} & 0 \leq \overline{S} - s_t \perp \overline{\sigma}_t \geq 0 \, ,&& \forall t \in \mathcal{T}\\
        \label{eq:kkt8} & 0 \leq p_t\uptxt{C} \perp \underline{\gamma}_t \geq 0 \, ,&& \forall t \in \mathcal{T}\\
        \label{eq:kkt9} & 0 \leq \overline{P}\uptxt{C} - p_t\uptxt{C} \perp \overline{\gamma}_t \geq 0 \, ,&& \forall t \in \mathcal{T}\\
        \label{eq:kkt10} & 0 \leq p_t\uptxt{D} \perp \underline{\delta}_t \geq 0 \, ,&& \forall t \in \mathcal{T}\\
        \label{eq:kkt11} & 0 \leq \overline{P}\uptxt{D} - p_t\uptxt{D} \perp \overline{\delta}_t \geq 0 \, ,&& \forall t \in \mathcal{T}.
    \end{align}
    \end{subequations}
We suppose that any optimal solution exhibits simultaneous charge and discharge, and prove a contradiction. We consider any $t \in \mathcal{T}$ such that $p_t\uptxt{C} > 0$ and $p_t\uptxt{D} > 0$. With \eqref{eq:kkt8} and \eqref{eq:kkt10}, $\underline{\gamma}_t = 0$ and $\underline{\delta}_t = 0$. We can then express $\lambda_t$ using \eqref{eq:kkt1}:
\begin{equation}
    \lambda_t = \eta\uptxt{D} ( \frac{- \overline{\delta}_t}{\Delta t} + C_t )
\end{equation}
We replace $\lambda_t$ in \eqref{eq:kkt2}, using that $\eta = \eta\uptxt{C} \eta\uptxt{D}$:
\begin{equation} \label{eq:contract}
    \Delta t \, C_t (1 - \eta) + \eta \overline{\delta}_t + \overline{\gamma}_t = 0,
\end{equation}
In the first case, $C_t > 0$ and $\eta < 1$, so $1 - \eta > 0$. We also have $\overline{\gamma}_t \geq 0$ and $\overline{\delta}_t \geq 0$ from \eqref{eq:kkt9} and \eqref{eq:kkt11}, so~\eqref{eq:contract} cannot hold. Therefore, an optimal solution cannot exhibit simultaneous charge and discharge.

Now suppose that $\eta = 1$, and again consider any $t \in \mathcal{T}$ such that $p_t\uptxt{C} > 0$ and $p_t\uptxt{D} > 0$. With \eqref{eq:contract}, we obtain that $\overline{\delta}_t =\overline{\gamma}_t=0$.
So the conditions \eqref{eq:kkt8} and \eqref{eq:kkt9} and the conditions \eqref{eq:kkt10} and \eqref{eq:kkt11} are valid for any $0 \leq p_t\uptxt{C} \leq \overline{P}\uptxt{C} $ and for any $0 \leq p_t\uptxt{D} \leq \overline{P}\uptxt{D} $ respectively. In particular, if, for $ s_t -  s_{t-1} \leq 0$, we take $p_t\uptxt{C} = 0 $ and $p_t\uptxt{D} = {-\eta\uptxt{D}(s_t -  s_{t-1})}/{\Delta t}$, all KKTs are satisfied, meaning that this solution is also optimal. If, for $ s_t -  s_{t-1} \geq 0$, we take $p_t\uptxt{C} = ({s_t -  s_{t-1}})/({\eta\uptxt{C} \Delta t}) $ and $p_t\uptxt{D} = 0$, all KKTs are satisfied, meaning that this solution is also optimal. Thus, there always exists an optimal solution that does not exhibit simultaneous charge and discharge.

Finally, we consider the case where $C_t\geq0$, $\forall t \in \mathcal{T}$, and $\eta < 1$. We focus on any $t \in \mathcal{T}$, such that $C_t=0$, and $p_t\uptxt{C} > 0$ and $p_t\uptxt{D} > 0$. Condition \eqref{eq:contract} reduces in the same way as before, so we can find another optimal solution without simultaneous charge and discharge using the same formulas for $p_t\uptxt{C}$ and $p_t\uptxt{D}$.

\section{Proof Corollary \ref{corol1}} \label{app:corol1}
Consider any $t \in \mathcal{T}$ such that there is simultaneous charge and discharge. We have  $1 - \eta > 0$ and  $\overline{\gamma}_t \geq 0$ from~\eqref{eq:kkt9}. Using \eqref{eq:contract}, $p_t\uptxt{C} > 0$ and $p_t\uptxt{D} > 0$ can only be satisfied if $C_t \leq 0$. We show in the proof of Proposition \ref{prop1} that if $C_t = 0$, we can construct another optimal solution without simultaneous charge and discharge at $t$. Therefore, simultaneous charge and discharge may only be unavoidable if $C_t < 0$.

\section{Proof Lemma \ref{lemma1}} \label{app:lemma1}
We introduce $\beta_t$ to represent the net charge, i.e., $\beta_t = s_t - \rho s_{t-1}$. We can rewrite \eqref{eq:update} at $t$ (or \eqref{eq:init} if $t=1$) as:
\begin{equation}
    \beta_t^* =  \Delta t \, \eta\uptxt{C} ( p_t\uptxt{C*} - \frac{1}{\eta\uptxt{D}\eta\uptxt{C}} p_t\uptxt{D*}).
\end{equation}
In the case where $\beta_t^* = \Delta t \, \eta\uptxt{C} \overline{P}\uptxt{C}$, we thus have
\begin{equation}
     p_t\uptxt{C*} = \frac{1}{\eta\uptxt{D}\eta\uptxt{C}} p_t\uptxt{D*} + \overline{P}\uptxt{C},
\end{equation}
which can only be feasible if $p_t\uptxt{D*} = 0$, meaning that there cannot be simultaneous charge and discharge. Similarly, we can prove that for $\beta_t^* = - \Delta t ({1}/{\eta\uptxt{D}}) \overline{P}\uptxt{D}$, $p_t\uptxt{C*} = 0$.

Now consider the case where $0 \leq \beta_t^* < \Delta t \, \eta\uptxt{C} \overline{P}\uptxt{C}$ and suppose that there is no simultaneous charge and discharge, i.e, $p_t\uptxt{D*}=0$ and $p_t\uptxt{C*} = {\beta_t^*}/({\Delta t \, \eta\uptxt{C}})$, with $p_t\uptxt{C*} < \overline{P}\uptxt{C}$.
We consider another solution $\mathbf{x}'$, which is the same, except for $0 < p_t\uptxt{D'} \leq \min\{\overline{P}\uptxt{D},\eta\uptxt{D}\eta\uptxt{C}(\overline{P}\uptxt{C}-p_t\uptxt{C*})\}$, and $p_t\uptxt{C'} = {\beta_t^*}/({\Delta t \, \eta\uptxt{C}})  +{p_t\uptxt{D'}}/({\eta\uptxt{D}\eta\uptxt{C}})  = p_t\uptxt{C*} + {p_t\uptxt{D'}}/({\eta\uptxt{D}\eta\uptxt{C}}) $. 
The existence of such a $p_t\uptxt{D'}$ is ensured since $\overline{P}\uptxt{D} > 0$ and $\overline{P}\uptxt{C}-p_t\uptxt{C*} > 0$.
The state of energy $s_t'$ is equal to $s_t^*$, 
so constraints \eqref{eq:init} to \eqref{eq:bounds_s} are satisfied by the new solution. 
Constraints \eqref{eq:bounds_pc} and \eqref{eq:bounds_pd} are satisfied at $t$ by definition of $p_t\uptxt{D'}$. This new solution is thus feasible.
We denote by $Z'$ the value of the objective function for the new solution, and we compare it to the value of the objective function for the initial solution, $Z^*$:
\begin{subequations}
    \allowdisplaybreaks
    \begin{align}
    Z' & = \sum_{t \in \mathcal{T}} \Delta t \, C_t (p_t\uptxt{D'}-p_t\uptxt{C'})   \\
    & = Z^* + \Delta t \, C_{t} (p_t\uptxt{D'}-p_t\uptxt{C'}) - \Delta t \, C_{t} (p_{t}\uptxt{D*}-p_{t}\uptxt{C*})\\
    & = Z^* + \Delta t \, C_{t} p_t\uptxt{D'} (1 - \frac{1}{\eta}),
    \end{align}
\end{subequations}
which is strictly greater than $Z^*$, 
since $C_{t} < 0$, $p_t\uptxt{D'} > 0$ and $\eta < 1$, indicating a contradiction.
Hence, in this case, any optimal solution will exhibit simultaneous charge and discharge.
Similarly, we can prove that for $- \Delta t ({1}/{\eta\uptxt{D}}) \overline{P}\uptxt{D} < \beta_t^* \leq 0$, any optimal solution will exhibit simultaneous charge and discharge.

\section{Proof Corollary \ref{corol2}} \label{app:corol2}
For $t \in \mathcal{T}$ such that $C_t < 0$, suppose that $s_t^*- \rho s_{t-1}^* = \Delta t \, \eta\uptxt{C} \overline{P}\uptxt{C}$. 
Using Condition 1, $s_t^*- \rho s_{t-1}^* > \overline{S} - \rho \underline{S}$, and constraint \eqref{eq:bounds_s} does not stand.
Now suppose that $s_t^*- \rho s_{t-1}^* = -\Delta t \, \eta\uptxt{D} \overline{P}\uptxt{D}$.
Using Condition 2, $s_t^*- \rho s_{t-1}^* < \underline{S} - \rho \overline{S}$, and constraint \eqref{eq:bounds_s} does not stand.
We thus have that $-\Delta t \, \eta\uptxt{D} \overline{P}\uptxt{D} < s_t^*- \rho s_{t-1}^* < \Delta t \, \eta\uptxt{C} \overline{P}\uptxt{C}$, and according to Lemma \ref{lemma1}, simultaneous charge and discharge is the best.

\section{Proof Theorem \ref{theorem1}} \label{app:thm1}
We suppose that any optimal solution $\mathbf{x}^*$ to the problem is such that there is no simultaneous charge and discharge, and we show a contradiction.
With Lemma \ref{lemma1}, we know that in each time period $t \in [\tau_1,\tau_2]$, we are either net charging the maximum or net discharging the maximum. 
There are thus two subsets of $[\tau_1,\tau_2]$, $\mathcal{T} \uptxt{C}$ and $\mathcal{T} \uptxt{D}$ such that $\mathcal{T} \uptxt{C} \cup \mathcal{T} \uptxt{D} = [\tau_1,\tau_2]$ and using \eqref{eq:s_calc} between $\tau_1$ and $\tau_2$, we get
\begin{align}
    \nonumber s_{\tau_2}^* = & \rho^{\tau_2-\tau_1+1} s_{\tau_1-1}^* + \Delta t  \eta \uptxt{C} \overline{P}\uptxt{C} \sum_{t \in \mathcal{T}\uptxt{C}} (\rho^{\tau_2-t}) \\
    & - \Delta t \frac{1}{\eta \uptxt{D}} \overline{P}\uptxt{D}\sum_{t \in \mathcal{T}\uptxt{D}} (\rho^{\tau_2-t})
\end{align}
Using $\tau_2-\tau_1+1=\overline{n}$, given Condition 2 we have that $s_{\tau_2}^* < \overline{S}$.

Suppose that $\mathcal{T} \uptxt{D} = \varnothing$. Then we have
\begin{subequations}
\begin{align}
    s_{\tau_2}^* & = \rho^{\overline{n}} s_{\tau_1-1}^* + \Delta t ( \eta \uptxt{C} \overline{P}\uptxt{C}  \sum_{t=1}^{\overline{n}} \rho^{\overline{n}-t}) \\
    & \geq \rho^{\overline{n}} \underline{S} + \Delta t ( \eta \uptxt{C} \overline{P}\uptxt{C} \sum_{t=1}^{\overline{n}} \rho^{\overline{n}-t}).
\end{align}
\end{subequations}
Using Condition~1, the constraint on the capacity of the storage system \eqref{eq:bounds_s} would be violated, so  $\mathcal{T} \uptxt{D} \neq \varnothing$.

We modify the optimal solution $\mathbf{x}^*$ to charge a very small quantity in a time period of discharge, and discharge a small quantity after $\tau_2$, when prices are positive again, so as to return on the same trajectory for the state of energy.
We identify the new solution with the superscript $'$.
Let's consider $\tau$ the last $t$ in $[\tau_1,\tau_2]$ such that $p_t\uptxt{D*} = \overline{P}\uptxt{D}$, and $\tau'$, the first $t > \tau_2$ such that $p_t\uptxt{D*} < \overline{P}\uptxt{D}$. We first suppose that $\tau'$ exists and then treat separately the case where $\tau'$ does not exist.
We introduce $\epsilon$ such that $0 < \epsilon \leq \min \{\overline{S}(1-\rho)+({1}/{\eta\uptxt{D}}) \overline{P}\uptxt{D}, ({\overline{S} - s_{\tau_2}^*})/{\rho^{\tau_2-\tau}}, \Delta t ({\overline{P}\uptxt{D}-p_{\tau'}\uptxt{D*}})/({\eta\uptxt{D} \rho^{\tau'-\tau}}) \}$, which is strictly greater than 0, so such an $\epsilon$ exists.
\begin{description}
    \item[For $t < \tau$:] We do not modify the solution.
    \item[For $t=\tau$:] We increase the charge, such that $p_\tau\uptxt{C'}={\epsilon}/({\Delta t \, \eta\uptxt{C}})$. The new state of energy is $s_\tau' = s_\tau^* + \epsilon$, which is still greater than $\underline{S}$.
    We also have that $s_\tau^* = \rho s_{\tau-1}^* - ({1}/{\eta\uptxt{D}}) \overline{P}\uptxt{D}$, and $s_{\tau-1}^* \leq \overline{S}$, so $s_\tau^* \leq \rho \overline{S}  - ({1}/{\eta\uptxt{D}}) \overline{P}\uptxt{D} $, and $s_\tau^*$ can be increased by a small quantity $\epsilon \leq \overline{S}(1-\rho)+({1}/{\eta\uptxt{D}}) \overline{P}\uptxt{D}$ without violating the upper bound.
    \item[For $\tau < t \leq \tau_2$:] We do not modify the solution.
    With \eqref{eq:s_calc} between $\tau$ and $t$ we get $s_t' = \rho^{t-\tau} s_\tau' + \Delta t $ $\sum_{k=\tau+1}^{t} (\rho^{t-k} \eta\uptxt{C} p_k\uptxt{C*}) = s_t^* + \rho^{t-\tau} \epsilon$, which is still greater than $\underline{S}$.
    Since in these time periods we are only charging, and with Assumption~\ref{ass3}, it is enough to prove that the upper bound on $s$ stands for the last time period, $\tau_2$. We have $s_{\tau_2}' = s_{\tau_2}^* + \rho^{\tau_2-\tau} \epsilon \leq \overline{S}$, by definition of $\epsilon$.
    \item[For $\tau_2 < t < \tau'$:] Using \eqref{eq:s_calc} between $\tau_2$ and $t$, we have that $s_t' = \rho^{t-\tau_2} s_{\tau_2}' - \Delta t \sum_{k=\tau_2+1}^{t} \rho^{t-k} ({1}/{\eta\uptxt{D}}) \overline{P}\uptxt{D}$, which is lower than $\overline{S}$, since $s_{\tau_2}' \leq \overline{S}$. We can further express $s_t' = s_t^* + \rho^{t-\tau} \epsilon$, so it is also greater than $\underline{S}$.
    \item[For $t = \tau'$:] We increase the discharge, $p_{\tau'}\uptxt{D'}=p_{\tau'}\uptxt{D*}+ \rho^{\tau'-\tau} \eta\uptxt{D} {\epsilon}/{\Delta t} \leq  \overline{P}\uptxt{D}$, by definition of $\epsilon$.  
    The new state of energy is $s_{\tau'}' = \rho s_{\tau'-1}'+ \Delta t ( \eta\uptxt{C} p_{\tau'}\uptxt{C'} - ({1}/{\eta\uptxt{D}}) p_{\tau'}\uptxt{D'} )= \rho s_{\tau'-1}^* + \rho^{\tau'-\tau} \epsilon + \Delta t ( \eta\uptxt{C} p_{\tau'}\uptxt{C*} - ({1}/{\eta\uptxt{D}}) p_{\tau'}\uptxt{D*} - {\rho^{\tau'-\tau} \epsilon}/{\Delta t} ) = s_{\tau'}^*$. 
    \item [For $t > \tau'$:] There is no change.
\end{description}
We have thereby proven that the new solution is feasible. We can calculate the new value of the objective function $Z'$, based on the previous value $Z^*$: $Z' = Z^* - C_\tau  {\epsilon}/{\eta\uptxt{C}} + C_{\tau'} \eta\uptxt{D} \rho^{\tau'-\tau} \epsilon  > Z^*$, since $C_\tau < 0$, $\epsilon >0$ and $C_{\tau'} \geq 0$. We found a feasible solution that increases the objective function, showing a contradiction.

We now consider the case where $\tau'$ does not exist, either because $\tau_2$ is the last time period in the horizon, or if $p_t\uptxt{D*} = \overline{P}\uptxt{D}$, $\forall t > \tau_2$. We can increase the charge at $\tau$ by a small quantity $\epsilon$, where $0 < \epsilon \leq \overline{S} - s_{\tau_2}^*$, without needing to decrease it later since the maximum level reached will be at $\tau_2$ and will be at most $\overline{S}$ by definition of $\epsilon$. 
The existence of $\epsilon$ is ensured since $s_{\tau_2}^* < \overline{S}$. 
In this case, it is clear that the objective function increases since the charge increases at $\tau$ and $C_\tau < 0$. Therefore, there is also a contradiction in this case.

\section{Proof Theorem \ref{theorem2}} \label{app:thm2}
With Corollary \ref{corol1}, it is sufficient to prove that simultaneous charge and discharge is not optimal in $[\tau_1,\tau_2]$. Consider an optimal solution $\mathbf{x}^*$ such that there exists at least one $t \in [\tau_1,\tau_2]$ such that $p_t\uptxt{C*} > 0$ and $p_t\uptxt{D*} > 0$.

We consider separately the case where $s_{t}^*<\overline{S},\forall t \in [\tau,\tau_2]$, and the case where $\exists \, t \in [\tau,\tau_2]$ such that $s_{t}^*=\overline{S}$,
and prove that in both cases, we can find a better solution $\mathbf{x}'$.

\textbf{First, in the case $s_{t}^*<\overline{S},\forall t \in [\tau,\tau_2]$,} we can increase the state of energy for some $\tau \in [\tau_1,\tau_2]$, and decrease it at a later time $\tau' > \tau_2$.
We choose $\tau$ such that $p_{\tau}\uptxt{D*} > 0$ and $p_{t}\uptxt{D*} = 0$, $\forall t \in [\tau+1,\tau_2]$.
The existence of $\tau$ is ensured as there exists at least one $t \in [\tau_1,\tau_2]$ such that $p_t\uptxt{D*} > 0$.
We choose $\tau'$ such that $p_{\tau'}\uptxt{D*} < \overline{P}\uptxt{D}$ and $p_{t}\uptxt{D*} = \overline{P}\uptxt{D}$ and $p_{t}\uptxt{C*} = 0$, $\forall t \in [\tau_2+1,\tau'-1]$.
We consider separately the case where such a $\tau'$ does not exist.
Both $\tau$ and $\tau'$ are represented in Figure \ref{fig:proof1}.

\ifarxiv
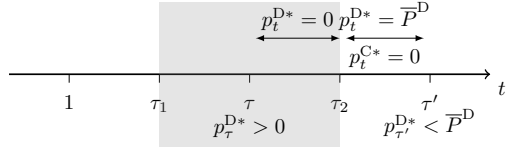
\begin{figure}[hb]
    \centering
    {\resizebox{.41\linewidth}{!}{\begin{tikzpicture}
    \draw[thick,->] (0,0) -- (8,0) node[below right] {$t$};

    \foreach \pos/\label/\shift in {1/1/0.5, 2.5/\tau_1/2.5, 4/\tau/2.5, 5.5/\tau_2/2.5, 7/\tau'/0} {
        \draw (\pos,0) -- ++ (0,-4pt) coordinate[label={[yshift=-2pt-\shift]below:$\label$}] {};
    }
    \fill[gray, fill opacity=0.2] (2.5,-1.2) rectangle (5.5,1.2);

    \draw[<->, >=latex] (5.6,0.6) -- node[above] {$p_t\uptxt{D*} = \overline{P}\uptxt{D}$} (6.9,0.6);
    \node[below] at (6.25,0.6) {$p_t\uptxt{C*} = 0$};
    \node[below] at (7,-0.45) {$p_{\tau'}\uptxt{D*} < \overline{P}\uptxt{D}$};

    \draw[<->, >=latex] (4.1,0.6) -- node[above] {$p_t\uptxt{D*} = 0$} (5.5,0.6);
    \node[below] at (4,-0.55) {$p_{\tau}\uptxt{D*} > 0$};
\end{tikzpicture}}}
    \caption{Different time periods for $s_{t}^*<\overline{S},\forall t \in [\tau,\tau_2]$. The interval with strictly negative prices is highlighted in gray.}
    \label{fig:proof1}
\end{figure}
\fi

For $t>\tau_2$, we know from Corollary \ref{corol1} that there is an optimal solution without simultaneous charge and discharge since prices are positive.
We decrease the discharge at $\tau$, such that $p_\tau\uptxt{D'}=p_\tau\uptxt{D*}-\epsilon$, with $\epsilon > 0$ and
\begin{equation} \label{eq:epsilon_1}
    \epsilon \leq \min \{ \frac{\overline{S}-s_{t}^*}{\Delta t \, \rho^{t-\tau}}\eta\uptxt{D},\forall t \in [\tau,\tau_2];p_\tau\uptxt{D*};\frac{\overline{P}\uptxt{D}-p_{\tau'}\uptxt{D*}}{\rho^{\tau'-\tau}}\}.
\end{equation}
Note that the existence of $\epsilon$ is ensured since all the elements considered are strictly positive.
We increase the discharge at $\tau'$, such that $p_{\tau'}\uptxt{D'}=p_{\tau'}\uptxt{D*}+\rho^{\tau'-\tau} \epsilon$.

We evaluate the feasibility of this new solution. The values of $p_t\uptxt{C}$ are not modified.
For $p_t\uptxt{D}$, we only have changes at $\tau$ and at $\tau'$, and the bounds in \eqref{eq:bounds_pd} are still respected, with \eqref{eq:epsilon_1}. We need to check that the bounds for the state of energy are respected for $\tau \leq t < \tau'$. It can be easily verified that for $t\geq \tau'$, $s_t^* = s_t'$.
\begin{description}
    \item[For $\tau \leq t \leq \tau_2$:] Using \eqref{eq:s_calc} between $\tau-1$ and $t$ and identifying $s_t^{*}$ we get:
        \begin{equation}
            s_t'= s_t^{*} + \Delta t \, \rho^{t-\tau} \frac{\epsilon}{\eta\uptxt{D}} ,
        \end{equation}
        which is greater than $\underline{S}$, since $s_t^{*}$ is.
        The upper bound on the state of energy is respected with \eqref{eq:epsilon_1}. 
    \item[For $\tau_2 < t < \tau'$:]  Using \eqref{eq:s_calc} between $\tau_2$ and $t$ we get:
                \begin{subequations}
                \allowdisplaybreaks
                    \begin{align}
                    s_t' & = \rho^{t-\tau_2} s_{\tau_2}' - \Delta t \sum_{k=\tau_2 + 1}^{t} \rho^{t-k} (\frac{1}{\eta\uptxt{D}} p_k\uptxt{D*} ) \\
                    & = s_{t}^{*} + \Delta t \, \rho^{t-\tau} \frac{\epsilon}{\eta\uptxt{D}} ,
                    \end{align}
                \end{subequations}
                which is both lower than $\overline{S}$, since $s_{\tau_2}' \leq \overline{S}$, and greater than $\underline{S}$, since $s_{t}^{*} \geq \underline{S}$.
\end{description}
We have thereby proven that in this case, and if $\tau'$ exists, the new solution is feasible. We now calculate the objective value $Z'$ and we compare it to the value of the objective function for the initial solution, $Z^*$:
    \begin{equation}
        Z' =  Z^* - \Delta t  \, C_{\tau} \, \epsilon + \Delta t \,  C_{\tau'} \, \epsilon \, \rho^{\tau'-\tau} > Z^*,
    \end{equation}
since $\epsilon > 0$, $C_{\tau} < 0$ and $C_{\tau'} \geq 0$.
The solution constructed is better, indicating a contradiction in this case.

Now we consider the case where $\tau'$ does not exist. This is the case if $p_t\uptxt{D*} = \overline{P}\uptxt{D}$, $\forall t > \tau_2$. We can still decrease discharge at $\tau$ by a small quantity $\epsilon$, such that $0< \epsilon \leq \min \{ ({\overline{S}-s_{t}^*})/({\Delta t \, \rho^{t-\tau}})\eta\uptxt{D},\forall t \in [\tau,\tau_2]; p_\tau\uptxt{D*} \}$. Until $\tau_2$, we can use the same proof to show that the storage bounds are not violated. If $\tau_2$ is not the last time period, for $t > \tau_2$, the state of energy is higher than in the initial solution, and it only decreases, due to discharging the maximum quantity in all the time periods, so we know that the bounds are not violated. The objective function for this modified solution is clearly higher since discharge decreases at $\tau$ and $C_\tau < 0$. Therefore, there is also a contradiction in the case $s_\tau^*<\overline{S}$ when $\tau'$ does not exist.

\textbf{We now look into the second case, i.e, $\exists \, t \in [\tau,\tau_2]$ such that $s_{t}^*=\overline{S}$.} We can decrease the state of energy at some $\tau'<\tau_1$ so that it can be increased at some $\tau \in [\tau_1,\tau_2]$. We choose $\tau$ such that $p_{\tau}\uptxt{D*} > 0$ and $p_{t}\uptxt{D*} = 0$, $\forall t \in [\tau_1,\tau-1]$. Its existence is ensured as there exists at least one $t \in [\tau_1,\tau_2]$ such that $p_t\uptxt{D*} > 0$. We choose $\tau'$ such that $p_{\tau'}\uptxt{D*} < \overline{P}\uptxt{D}$, and $p_{t}\uptxt{D*} = \overline{P}\uptxt{D}$ and $p_{t}\uptxt{C*} = 0$, $\forall t \in [\tau'+1,\tau_1-1]$. All are represented in Figure \ref{fig:proof2}.
\ifarxiv
\begin{figure}[hb]
    \centering
    {\resizebox{.41\linewidth}{!}{\begin{tikzpicture}
    \draw[thick,->] (0,0) -- (8,0) node[below right] {$t$};

    \foreach \pos/\label/\shift in {1/1/0.5, 2.5/\tau'/0, 4/\tau_1/2.5, 5.5/\tau/2.5, 7/\tau_2/2.5} {
        \draw (\pos,0) -- ++ (0,-4pt) coordinate[label={[yshift=-2pt-\shift]below:$\label$}] {};
    }
    \fill[gray, fill opacity=0.2] (4,-1.2) rectangle (7,1.2);

    \draw[<->, >=latex] (2.6,0.6) -- node[above] {$p_t\uptxt{D*} = \overline{P}\uptxt{D}$} (3.9,0.6);
    \node[below] at (3.25,0.6) {$p_t\uptxt{C*} = 0$};
    \node[below] at (2.5,-0.45) {$p_{\tau'}\uptxt{D*} < \overline{P}\uptxt{D}$};

    \draw[<->, >=latex] (4,0.6) -- node[above] {$p_t\uptxt{D*} = 0$} (5.4,0.6);
    \node[below] at (5.5,-0.55) {$p_{\tau}\uptxt{D*} > 0$};
\end{tikzpicture}}}
    \caption{Different time periods when $\exists \, t \in [\tau,\tau_2]$ such that $s_{t}^*=\overline{S}$. The interval with strictly negative prices is in gray.}
    \label{fig:proof2}
\end{figure}
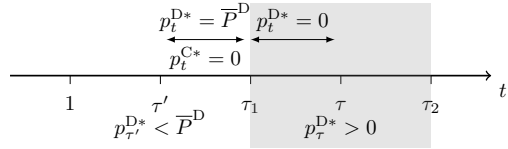
\fi

We prove the existence of $\tau'$ by contradiction. Suppose that for all $t<\tau_1$, $p_t\uptxt{D*} = \overline{P}\uptxt{D}$.
Then, with \eqref{eq:s_calc} between $0$ and $\tau_1-1$
\begin{align}
    s_{\tau_1-1}^* & = \rho^{\tau_1-1} S\uptxt{init} - \Delta t (\sum_{t=1}^{\tau_1-1} \rho^{\tau_1-1-t}) \frac{1}{\eta\uptxt{D}} \overline{P}\uptxt{D},
\end{align}
and using \eqref{eq:s_calc} between $\tau_1-1$ and $\tau_2$
\begin{align}
    \nonumber s_{\tau_2}^* & = \rho^{\tau_2-\tau_1+1} s_{\tau_1-1}^* \\
    & \quad + \Delta t \sum_{t=\tau_1}^{\tau_2} \rho^{\tau_2-t} ( \eta\uptxt{C} p_t\uptxt{C*} - \frac{1}{\eta\uptxt{D}} p_t\uptxt{D*} ).
\end{align}
We know that $\sum_{t=\tau_1}^{\tau_2} \rho^{\tau_2-t} ( \eta\uptxt{C} p_t\uptxt{C*} - ({1}/{\eta\uptxt{D}}) p_t\uptxt{D*} ) < \sum_{t=\tau_1}^{\tau_2} \rho^{\tau_2-t} ( \,\eta\uptxt{C} \overline{P}\uptxt{C} )$, since there is at least one $t \in [\tau_1,\tau_2]$ such that $p_t\uptxt{D*} > 0$.
Using the second condition, we get the result $s_{\tau_2}^* < \overline{S}$, indicating a contradiction. The existence of $\tau'$ is thus guaranteed.

We decrease the discharge at $\tau$, such that $p_\tau\uptxt{D'}=p_\tau\uptxt{D*}-\epsilon$, with 
\begin{equation} \label{eq:epsilon_2}
    0 < \epsilon \leq \min \{ \frac{s_{\tau_1-1}^*-\underline{S}}{\Delta t\, \rho^{\tau_1-\tau'-1}}\eta\uptxt{D},\frac{p_\tau\uptxt{D*}}{\rho^{\tau-\tau'}},\overline{P}\uptxt{D}-p_{\tau'}\uptxt{D*} \}.
\end{equation}
To be able to do this, we increase the discharge at $\tau'$, such that $p_{\tau'}\uptxt{D'}=p_{\tau'}\uptxt{D*}+\rho^{\tau-\tau'}\epsilon$.
Note that the existence of $\epsilon$ is ensured since all the elements are strictly positive. For the first element, the strict positivity comes from the fact that $s_{\tau-1}^*>\underline{S}$, otherwise $s_{\tau_2}^*=\overline{S}$ could not be reached, due to the second condition.

We evaluate the feasibility of this new solution. The values of $p_t\uptxt{C}$ are not modified.
For $p_t\uptxt{D}$, we only have changes at $\tau$ and at $\tau'$, and the bounds in \eqref{eq:bounds_pd} are still respected, with \eqref{eq:epsilon_2}. We need to check that the bounds for the state of energy are respected for $\tau' \leq t < \tau$. It can be easily verified that for $t\geq \tau$, $s_t^* = s_t'$.
\begin{description}
    \item[For $\tau' \leq t < \tau_1$:] Using \eqref{eq:s_calc} between $\tau'-1$ and $t$, we get:
                \begin{align}
                    s_t' = s_t^{*} - \Delta t \rho^{t-\tau'} \frac{1}{\eta\uptxt{D}} \epsilon,
                \end{align}
                which is lower than $\overline{S}$, since $s_t^{*}$ is.
                It is enough to check that the lower bound on the state of energy is respected at $\tau_1-1$ since there is only discharge between $\tau'$ and $\tau_1-1$. It is the case since $\epsilon \leq  ({s_{\tau_1-1}^*-\underline{S}})/({\Delta t \, \rho^{\tau_1-\tau'-1}})\eta\uptxt{D}$.
    \item[For $\tau_1 \leq t < \tau$:]  Using \eqref{eq:s_calc} between $\tau_1-1$ and $t$, we get:
                \begin{subequations}
                \allowdisplaybreaks
                    \begin{align}
                    s_t' & = \rho^{t-\tau_1+1} s_{\tau_1-1}' + \Delta t \sum_{k=\tau_1}^{t} \rho^{t-k} \eta\uptxt{C} p_k\uptxt{C*}  \\
                    & = s_{t}^{*} - \Delta t \, \rho^{t-\tau'} \frac{1}{\eta\uptxt{D}} \epsilon,
                    \end{align}
                \end{subequations}
                which is both greater than $\underline{S}$, since $s_{\tau_1-1}' \geq \underline{S}$, and lower than $\overline{S}$, since $s_{t}^{*} \leq \overline{S}$.
\end{description}
We have thereby proven that in the case when $\exists \, t \in [\tau,\tau_2]$ such that $s_{t}^*=\overline{S}$, the new solution is feasible. We can calculate the objective value $Z'$ and compare it to $Z^*$:
    \begin{equation}
        Z' =  Z^* + \Delta t  C_{\tau'} \epsilon - \Delta t  C_{\tau} \epsilon \rho^{\tau-\tau'} > Z^*,
    \end{equation}
since $\epsilon > 0$, $C_{\tau} < 0$ and $C_{\tau'} \geq 0$.
The solution constructed is better, indicating a contradiction in this case.
    
\section{Proof Theorem \ref{theorem3}} \label{app:thm3}

Suppose that the conditions of Theorem \ref{theorem3} hold and consider $\mathbf{x}^*$, an optimal solution to \eqref{pb:opt} with $p_\tau\uptxt{D*} > 0$ and $p_\tau\uptxt{C*} > 0$, where $\tau \in \mathcal{T}$. With Corollary \ref{corol1}, $C_\tau < 0$. 

We consider separately the two possible situations:
\begin{enumerate}
    \item $s_t^* < \overline{S}$, $\forall t \geq \tau$
    \item $\exists \, t \geq \tau$ such that $s_t^* = \overline{S}$.
\end{enumerate}

\begin{description}
    \item[Case 1:] In the first case we can decrease the discharge at $\tau$ by a little quantity $\epsilon > 0$, without violating the upper bound on the state of energy for the rest of the horizon. The objective function is increased, since we subtract $C_\tau \, \epsilon$, which is strictly negative. This contradicts the fact that $\mathbf{x}^*$ is optimal.

    \item[Case 2:] We introduce $t\uptxt{max}$, the first $t \geq \tau$ such that $s_{t\uptxt{max}}^* = \overline{S}$. There are two possible sub-cases:
    \begin{itemize}
        \item [2.a:] $t\uptxt{max}\in\mathcal{T^+}$
        \item [2.b:] $t\uptxt{max}\in\mathcal{T^-}$
    \end{itemize}

        \item[Case 2.a:] 
        If $t\uptxt{max}\in\mathcal{T^+}$, we can decrease the discharge at $\tau$ by a small quantity $\epsilon > 0$ and decrease the charge at $t\uptxt{max}\in\mathcal{T^+}$ by $\rho^{t\uptxt{max}-\tau} {\epsilon}/({\eta\uptxt{C}\eta\uptxt{D}})$. We know that the upper bound on the storage level will not be violated between $\tau$ and $t\uptxt{max}$ if $\epsilon$ is sufficiently small. The new value of the objective function is:
            \begin{equation}
                Z'= Z^* - C_{\tau} \epsilon + \rho^{t\uptxt{max}-\tau} C_{t\uptxt{max}}  \frac{\epsilon}{\eta\uptxt{C}\eta\uptxt{D}},
            \end{equation}
            so $Z'> Z^*$, since $C_{t\uptxt{max}} \geq 0$ and $ C_{\tau}<0$. This new solution is feasible and better, showing a contradiction. In this case, simultaneous charge and discharge cannot be optimal.
        \item[Case 2.b:]
        If $t\uptxt{max}\in\mathcal{T^-}$, we introduce $\nu \in \mathcal{J}$ such that $t\uptxt{max}\in \mathcal{T}^-_\nu$.
        We can prove that there exists $\tau' \in \mathcal{T^+}$ such that $\tau' < t\uptxt{max}$ and $p_{\tau'}\uptxt{D*} < \overline{P}\uptxt{D}$.
        Suppose that $\tau'$ does not exist, which means that $p_{t}\uptxt{D*} = \overline{P}\uptxt{D}$, $\forall t \in \mathcal{T^+}$, $t<t\uptxt{max}$. 
        Using~\eqref{eq:s_calc} between $0$ and $t\uptxt{max}$, expressing and isolating the terms related to $n_1$ and $p_1$, we can find an upper bound to $s_{t\uptxt{max}}^*$ involving $\hat{S}_1$. By induction until $\nu-1$, we get:
        \begin{equation} \label{pb:23}
                s_{t\uptxt{max}}^* <  \rho^m S + \Delta t \eta\uptxt{C} \overline{P}\uptxt{C} \frac{1-\rho^{m}}{1-\rho},
        \end{equation}
        where the strict inequality comes from $p\uptxt{D*}_\tau>0$, $m$ is such that $t\uptxt{max} = p_1 + n_1 +...+p_\nu + m$, and $S = (\rho^{p_\nu} \hat{S}_{\nu-1} - \Delta t ({1}/{\eta\uptxt{D}}) \overline{P}\uptxt{D} (\sum_{t=1}^{p_\nu}\rho^{p_\nu-t}) )$.

        Note that $S \leq \overline{S}$, since it corresponds to the state of energy obtained by discharging the maximum quantity during $p_\nu$ time periods, starting from $\hat{S}_{\nu-1} \leq \overline{S}$. For $\rho = 1$, it is immediate that $s_{t\uptxt{max}}^* < \hat{S}_\nu$, therefore $s_{t\uptxt{max}}^* < \overline{S}$ with Condition 1. 
        We consider that $\rho < 1$. From Assumption~\ref{ass3}, $S (1-\rho) \leq \Delta t \eta\uptxt{C} \overline{P}\uptxt{C}$. Multiplying both sides by $({\rho^m-\rho^{n_\nu}})/({1-\rho})$ (with $m\leq n_\nu$) and rearranging, we get
        \begin{equation}
            \rho^m S + \Delta t \eta\uptxt{C} \overline{P}\uptxt{C} \frac{1-\rho^{m}}{1-\rho} \leq \hat{S}_\nu.
        \end{equation}
        Using Condition 1, we get that $s_{t\uptxt{max}}^* < \overline{S}$, which is a contradiction. Therefore, $\tau'$ must exist.
        If several time periods correspond to this description, we take $\tau'$ to be the largest of those, so that we have $p_{t}\uptxt{D*} = \overline{P}\uptxt{D}$, $\forall t \in\mathcal{T^+}$, with $\tau' < t < t\uptxt{max}$.
        Here, we identify two possible situations:
        \begin{itemize}
            \item [2.b.i:] $\tau' < \tau$
            \item [2.b.ii:] $\tau' > \tau$
        \end{itemize}
        \item[Case 2.b.i:]
        We first show that $\underline{S}$ is not reached between $\tau'$ and $\tau$.
        We suppose that there exists $t\uptxt{min} \in \mathcal{T}$, with $\tau' < t\uptxt{min} < \tau$, such that $s_{t\uptxt{min}}^*=\underline{S}$, and show a contradiction. 
        Note that $t\uptxt{min}+1$ has to be in $\mathcal{T^-}$, since $p_{t}\uptxt{D*} = \overline{P}\uptxt{D}$, $\forall t \in\mathcal{T^+}$, and the minimum is already reached at $t\uptxt{min}$. Let's consider that $t\uptxt{min}+1 \in \mathcal{T}_\sigma^-$. The number of time periods left in $\mathcal{T}_\sigma^-$, including $t\uptxt{min}+1$, is $n'_\sigma$, and the number of periods already passed is $m'_\sigma$, such that $n_\sigma = m'_\sigma + n'_\sigma$.
        We express \eqref{eq:s_calc} between $t\uptxt{min}$ and $t\uptxt{max}$, and get an upper bound using $\eta\uptxt{C} p_t\uptxt{C*} - ({1}/{\eta\uptxt{D}}) p_t\uptxt{D*}$ with $\overline{P}\uptxt{C}$. Isolating the terms related to $n'_\sigma$, factoring $\rho^{m'_\sigma}$ in, and using $n'_\sigma \leq n_\sigma$, we get:
        \begin{align}
            \nonumber s_{t\uptxt{max}}^{*} & \leq \rho^{t\uptxt{max}-t\uptxt{min}-n_\sigma}  ( \rho^{n_\sigma} \underline{S}  \\
            \nonumber &  + \Delta t (\sum_{t=1}^{n_\sigma} \rho^{n_\sigma-t}) \eta\uptxt{C} \overline{P}\uptxt{C} ) \\
            & \quad + \Delta t\sum_{t=t\uptxt{min}+n'_\sigma+1}^{t\uptxt{max}} \rho^{t\uptxt{max}-t} ( \eta\uptxt{C} p_t\uptxt{C*} - \frac{1}{\eta\uptxt{D}} p_t\uptxt{D*} ).
        \end{align}
        Using that $\underline{S} \leq \max \{ \rho^{p_\sigma} \hat{S}_{\sigma-1} - \Delta t ({1}/{\eta\uptxt{D}}) \overline{P}\uptxt{D} (\sum_{t=1}^{p_\sigma} \rho^{p_\sigma-t}), \underline{S}\}$, and that $t\uptxt{max}-t\uptxt{min}-n'_\sigma-...-p_{\nu-1}-n_{\nu-1}-p_\nu=m$, we get:
            \begin{align} \label{pb:smin_end}
                \nonumber s_{t\uptxt{max}}^* & <  \rho^m (\rho^{p_\nu} \hat{S}_{\nu-1} + \Delta t \frac{1}{\eta\uptxt{D}} \overline{P}\uptxt{D} (\sum_{t=1}^{p_\nu} \rho^{p_\nu-t}) )\\
                & \quad + \Delta t \eta\uptxt{C} \overline{P}\uptxt{C} (\sum_{t=1}^{m} \rho^{m-t}),
            \end{align}
        so $s_{t\uptxt{max}}^* < \hat{S}_\nu \leq \overline{S}$, which is a contradiction. Therefore, $\underline{S}$ cannot be reached after $\tau'$. We can thus increase the discharge at $\tau'$ by a small quantity $\epsilon > 0$ without violating the lower bound on the storage level before~$\tau$. We decrease the discharge at $\tau$ by $\rho^{\tau-\tau'} \epsilon$, returning on the same trajectory, ensuring that the upper bound on the storage level is not violated at $\tau$ or after. The new value of the objective function is:
        \begin{equation}
            Z'= Z^* + C_{\tau'} \epsilon - \rho^{\tau-\tau'} C_{\tau} \epsilon,
        \end{equation}
        so $Z'> Z^*$, since $C_{\tau'} \geq 0$ and $ C_{\tau}<0$. This new solution is feasible and better, showing a contradiction. In this case, simultaneous charge and discharge cannot be optimal.
    
        \item[Case 2.b.ii:]
        If $\tau' > \tau$, since the maximum level is not reached before $t\uptxt{max}$, we know we can increase the level between $\tau$ and $\tau'$ by a small quantity $\Delta t \, ({1}/{\eta\uptxt{D}}) \epsilon > 0$ without violating the upper bound.
        We decrease the discharge at $\tau$ by $\epsilon$ and increase the discharge at $\tau'$ by $\rho^{\tau'-\tau} \epsilon$ to return on the same trajectory for the state of charge before $t\uptxt{max}$.
        The new value of the objective function is:
            \begin{equation}
                Z'= Z^* -  C_{\tau} \epsilon + \rho^{\tau'-\tau} C_{\tau'} \epsilon ,
            \end{equation}
            so $Z'> Z^*$, since $C_{\tau'} \geq 0$ and $ C_{\tau}<0$.
        The new solution is feasible and better, showing a contradiction, so simultaneous charge and discharge cannot be optimal, which completes the proof.
\end{description}
\end{appendices}

\bibliographystyle{abbrvnat}
\bibliography{ref}

\end{document}